\input amstex\documentstyle{amsppt}  
\pagewidth{12.5cm}\pageheight{19cm}\magnification\magstep1
\topmatter
\title Unipotent elements in small characteristic, II\endtitle
\author G. Lusztig\endauthor
\address{Department of Mathematics, M.I.T., Cambridge, MA 02139}\endaddress
\thanks{Supported in part by the National Science Foundation.}\endthanks
\endtopmatter   
\document

\define\dy{\dot y}

\define\mpb{\medpagebreak}

\define\da{\dagger}

\define\Lie{\text{\rm Lie }}

\define\frl{\forall}
\define\pe{\perp}
\define\si{\sim}

\define\sqc{\sqcup}

\define\qua{\quad}

\define\dx{\dot x}

\define\op{\oplus}

\define\part{\partial}
\define\em{\emptyset}
\define\imp{\implies}
\define\ra{\rangle}
\define\n{\notin}
\define\iy{\infty}
\define\m{\mapsto}
\define\do{\dots}
\define\la{\langle}

\define\sub{\subset}    

\define\T{\times}
\define\ti{\tilde}
\define\nl{\newline}
\redefine\i{^{-1}}

\define\Ad{\text{\rm Ad}}
\define\Hom{\text{\rm Hom}}
\define\End{\text{\rm End}}

\define\Ker{\text{\rm Ker}}

\redefine\d{\delta}
\define\e{\epsilon}
\define\et{\eta}

\redefine\o{\omega}

\define\ph{\phi}
\define\ps{\psi}
\define\r{\rho}

\redefine\l{\lambda}

\define\x{\xi}

\define\Ph{\Phi}
\define\Ps{\Psi}

\define\kk{\bold k}

\define\ZZ{\bold Z}

\define\cc{\Cal C}

\define\ce{\Cal E}
\define\cf{\Cal F}

\define\cj{\Cal J}
\define\ck{\Cal K}

\define\cm{\Cal M}
\define\cn{\Cal N}

\define\car{\Cal R}

\define\cv{\Cal V}
\define\cw{\Cal W}

\define\cx{\Cal X}
\define\cy{\Cal Y}

\define\fg{\frak g}

\define\fB{\frak B}

\define\fP{\frak P}

\define\fX{\frak X}

\define\ta{\ti a}

\define\tN{\ti N}

\define\tY{\ti Y}

\define\sps{\supset}

\define\tcf{\ti{\cf}}

\define\na{\nabla}
\define\vtr{\vartriangle}
\define\tcm{\ti\cm}
\define\vda{\vdash}
\define\bacc{\bar{\cc}}
\define\bV{\bar V}
\define\bQ{\bar Q}
\define\bR{\bar R}
\define\Rad{\text{\rm Rad}}
\define\Iso{\text{\rm Iso}}
\define\bN{\bar N}
\define\gr{\text{\rm gr}}
\define\KO{K}
\define\NU{L1}
\define\LU{L2}
\define\WA{W}
\head Introduction\endhead
\subhead 0.1\endsubhead
Let $\kk$ be an algebraically closed field of characteristic exponent $p\ge1$. This paper 
is a study of unipotent elements in a special orthogonal group $SO_Q$ (where $Q$ is a
nondegenerate quadratic form on a finite dimensional $\kk$-vector space $V$) with emphasis
on the case where $p=2$. We develop some of the proposals in \cite{\LU} which try to extend
the Dynkin-Kostant theory \cite{\KO} of unipotent elements in the case $p=1$ to the case 
$p>1$. 

Namely we show that to any unipotent element $u\in G$ one can associate canonically a 
filtration of $V$ whose is stabilizer in $SO_Q$ is a parabolic subgroup containing $u$ in 
the unipotent radical and such that this parabolic is of the same type as a parabolic
attached in the Dynkin-Kostant theory to a unipotent element in a special orthogonal group
with $p=1$. This allows us to partition the unipotent variety of $SO_Q$ into
pieces which are both smooth varieties and unions of (possibly several) unipotent conjugacy
classes. We show that (for a split rational structure over a finite field with $q$ 
elements) each piece has a number of rational points given by a polynomial in $q$ with 
integer coefficients independent of $q$ and $p$. (This kind of property was stated without
proof in \cite{\NU}.) The similar statement was established for groups of type $C$ in 
\cite{\LU}; for  an exceptional type this is is easily established using the available 
tables since in that case each piece contains at most two unipotent conjugacy classes. In 
some sense the behaviour of the pieces is more complicated for orthogonal groups than for 
other almost simple groups. For example this is the only case where the number of unipotent
conjugacy classes in a given piece (in bad characteristic) is not necessarily a power of 
$2$. 

{\it Notation.} We set $\ZZ'=2\ZZ+1,\ZZ''=2\ZZ$. For $a\in\ZZ$ we set
$\ZZ'_{\ge a}=\{n\in\ZZ';n\ge a\}$, $\ZZ''_{\ge a}=\{n\in\ZZ'';n\ge a\}$. The cardinal of a
finite set $X$ is denoted by $|X|$; if $X$ is an infinite set we write $|X|=\iy$.

\head Contents\endhead
1. Unipotent elements in orthogonal groups.

2. The case $p=2$.

3. The case $p\ne2$.

4. On unipotent conjugacy classes in $SO_Q$ ($p=2$).

\head 1. Unipotent elements in orthogonal groups\endhead
\subhead 1.1\endsubhead
Let $\kk$ be a finite or algebraically closed field of characteristic exponent $p\ge1$. Let
$\cc$ be the category whose objects are $\kk$-vector spaces of finite dimension; the
morphisms are linear maps. Let $\cc'$ be the subcategory of $\cv$ consisting of all 
$V\in\cv$ such that $\dim V\in\ZZ'$. Let $\cc''$ be the subcategory of $\cv$ consisting of
all $V\in\cv$ such that $\dim V\in\ZZ''$. Let $\bacc$ be the category whose objects are 
$\ZZ$-graded $\kk$-vector spaces $\bV=\op_{a\in\ZZ}\bV^a$ such that $\dim\bV<\iy$; the 
morphisms are linear maps respecting the grading.

Let $V\in\cc$. For $N\in\End V$ nilpotent and $i\ge1$ let

$E_i^N=\ker N^i/(N(\ker N^{i+1})+\ker N^{i-1})$.
\nl
Then $c_i^N:=\dim E_i^N$ is the number of Jordan blocks of size $i$ of $N$. Let $e=e_N$ be 
the smallest integer $\ge0$ such that $N^e=0$.

For any symmetric or symplectic bilinear form $\la,\ra:V\T V@>>>\kk$ let
$\Rad_{\la,\ra}=\{x\in V;\la x,V\ra=0\}$ be the radical of $\la,\ra$.

In the remainder of this paper we fix $V\in\cc$ and a quadratic form $Q:V@>>>\kk$ with 
associated symmetric bilinear form $\la,\ra:V\T V@>>>\kk$; we have
$\la x,y\ra=Q(x+y)-Q(x)-Q(y)$ for $x,y\in V$. Let $D=\dim V$. Let $R=\Rad_{\la,\ra}$. We 
assume that $Q$ is {\it nondegenerate} that is, $Q|_R:R@>>>\kk$ is injective. (In this case
we must have $R=0$ unless $V\in\cc',p=2$ when $\dim R=1$.)

In the case where $D\in\ZZ''$ let $\cj_Q$ be the set of $D/2$-dimensional subspaces of $V$
on which $Q$ is $0$. In the case where $|\kk|<\iy,D\in\ZZ''_{\ge2}$ we set $\et_Q=1$ if 
$\cj_Q\ne\em$ and $\et_Q=-1$ if $\cj_Q=\em$.

For any subspace $W$ of $V$ we set $W^\pe=\{x\in V;\la x,W\ra=0\}$. 

Let $O_Q=\{T\in GL(V);Q(Tx)=Q(x)\qua\frl x\in V\}$ be the orthogonal group of $Q$. In the 
case where $|\kk|=\iy$ we denote by $SO_Q$ the identity component of $O_Q$; in the case 
where $|\kk|<\iy$ we set $SO_Q=O_Q\cap SO_{\ti Q}$ where $\ti Q$ is the quadratic form 
obtained from $Q$ by extension of scalars to an algebraic closure of $\kk$. Let
$$\align\tcm_Q&=\{N\in\End(V);1+N\text{ unipotent in }O_Q\}\\&
=\{N\in\End(V);N\text{ nilpotent },Q(Nx)=-\la x,Nx\ra\text{ for all }x\in V\}.\endalign$$
Note that if $N\in\tcm_Q$ and $x,y\in V$ then
$$\la x,Ny\ra+\la Nx,y\ra+\la Nx,Ny\ra=0.$$
In particular, for $N\in\tcm_Q$ we have $\la x,Ny\ra=\la N^\da x,y\ra$ for $x,y\in V$, 
where $N^\da:=(1+N)\i-1=-N+N^2-N^3+\do\in\tcm_Q$. 

Let $N\in\tcm_Q$. We have
$$NV\cap R=0.\tag a$$
If $x\in V$ and $Nx\in R$ then $Q(Nx)=-\la x,Nx\ra=0$; since $Q:R@>>>\kk$ is injective we 
see that $Nx=0$. Thus (a) holds. 

From the definitions we have for $i\ge1$:

(b) $(\ker N^i)^\pe=N^iV+R$, $(N^iV)^\pe=\ker N^i$.
\nl
In particular, $NR\sub R$. Since $\dim R\le1$, it follows that $NR=0$.

\subhead 1.2\endsubhead
In this subsection we assume that $p>1$ and that $|\kk|=q<\iy$.
For any $m\in ZZ'_{\ge1}$ we set

$P_m=q^{1/4(m^2-2m+1)}(q^2-1)(q^4-1)\do(q^{m-1}-1)$.
\nl
For any $m\in\ZZ''_{\ge2}$ and $\d\in\{1,-1\}$ we set

$P_m^\d=q^{1/4(m^2-2m)}(q^2-1)(q^4-1)\do(q^{m-2}-1)(q^{m/2}-\d)$.
\nl 
Let $\cv$ be a $\kk$-vector space with a nondegenerate quadratic form $u:\cv@>>>\kk$. Let 
$(,):\cv\T\cv@>>>\kk$ be the associate symmetric bilinear form. Let $\car=\Rad_{(,)}$. Let
$s=\dim\cv$. Let $k\in[0,s]$. Let $S_k\cv$ be the set of $k$-dimensional subspaces $\cw$ of
$\cv$ such that $u|_\cw$ is nondegenerate. If $k\in\ZZ''_{\ge2}$ we have a partition
$S_k\cv=\sqc_{\d\in\{1,-1\}}S_k\cv^\d$ where $S_k\cv^\d=\{\cw\in S_k\cv;\et_{u_\cw}=\d\}$.
If $s\in\ZZ''$, $k\in\ZZ'$ we set $N_{s,k}^{\e,*}=|S_k\cv|$. If $s\in\ZZ'$, $k\in\ZZ'$, we
set $N_{s,k}^{*,*}=|S_k\cv|$. If $s\in\ZZ''$, $k\in\ZZ''_{\ge2}$ and $\d\in\{1,-1\}$, we 
set $N_{s,k}^{\e,\d}=|S_k^\d\cv|$. If $s\in\ZZ'$, $k\in\ZZ''_{\ge2}$ and $\d\in\{1,-1\}$, 
we set $N_{s,k}^{*,\d}=|S_k^\d\cv|$.

(i) {\it If $s\in\ZZ''_{\ge2}$, $k\in\ZZ'$, $\e=\et_u$ then 
$N_{s,k}^{\e,*}=P_s^\e P_k\i P_{t-k}\i$.}
\nl
If $p\ne2$ this is obvious. Assume now that $p=2$. Let $\car'$ be the radical of
$(,)|_\cw$ for $\cw\in S_k\cv$. Then $\car'$ is one of the $q^{(s-2)/2}(q^{s/2}-\e)$ lines
in $\cv$ on which $u$ is not identically $0$. Also $\car'{}^\pe$ is a hyperplane in $\cv$ 
and $(,)$ induces a nondegenerate symplectic form on $\car'{}^\pe/\car'$. For each such 
$\car'$, the number of $\cw\in S_k\cv$ such that $\car'\sub\cw\sub\car'{}^\pe$ is the
number of $(k-1)$-dimensional subspaces of $\car'{}^\pe/\car'$ on which $(,)$ is 
nondegenerate hence it is $P_{s-1}P_k\i P_{s-k}\i$. Hence 
$N_{s,k}^{\e,*}=q^{(s-2)/2}(q^{s/2}-\e)P_{s-1}P_k\i P_{s-k}\i$, as required.

(ii) {\it If $s\in\ZZ''_{\ge2}$, $k\in\ZZ''_{\ge2}$, $k<s$, $\e=\et_u$, $\d\in\{1,-1\}$, 
then 

$N_{s,k}^{\e,\d}=2\i P_s^\e(P_k^\d)\i(P_{s-k}^{\e\d})\i$.}
\nl
This is clear: the orthogonal group of $u$ acts transitively on $S_k^\d\cv$.

(iii) {\it If $s\in\ZZ'$, $k\in\ZZ'$, $k<s$, then 
$N_{s,k}^{*,*}=2\i\sum_{\d\in\{1,-1\}}P_sP_k\i(P_{s-k}^\d)\i$.}
\nl
If $p\ne2$ this is obvious. Assume now that $p=2$. We set $s=2t+1,k=2a+1$. Let 
$\car'=\Rad_{(,)|_\cw}$ for $\cw\in S_k\cv$. Then $\car'$ is one of the $q^{2t}$ lines in 
$\cv$ on which $u$ is not identically $0$. The number of $\cw\in S_k\cv$ such that the 
corresponding $\car'$ is equal to $\car$ is the number of $2a$-dimensional subspaces of 
$\cv/\car$ on which $(,)$ is nondegenerate that is $P_sP_k\i P_{s-k+1}\i$. Let $N'$ be the
number of $\cw\in S_k\cv$ such 
that the corresponding $\car'$ is not equal to $\car$. Then $N'=(q^{s-1}-1)N'_0$ where 
$N'_0$ is the number of $\cw$ with prescribed $\car'\ne\car$. To such $\cw$ we associate 
the subspace $\cw'=\cw+\car$ (of dimension $k+1$) of $\car'{}^\pe$ (of dimension $s-1$). 
Note that the number of possible such $\cw'$ is the number of $k-1$-dimensional subspaces
of $\car'{}^\pe/(\car+\car')$ on which $(,)$ is nondegenerate that is 
$P_{s-2}P_k\i P_{s-k-1}\i$. For any $\cw'$ as above the number of $\cw$ such that 
$\car'\sub\cw\sub\cw'$ and $\cw\op\car=\cw'$ is the number of subspaces $\cx\sub\cw'/\car'$
such that $\cx\op\car=\cw'/\car'$ that is $q^{k-1}$. Thus,
$$\align&N'=(q^{s-1}-1)q^{k-1}P_{s-2}P_k\i P_{s-k-1}\i\\&
=q^{2a}q^{-t^2+(t-1)^2}q^{(t-a)^2-(t-a-1)^2}(q^{2t-2a}-1)P_sP_k\i P_{s-k+1}\i.\endalign$$
We see that
$$\align N_{s,k}^{*,*}&=P_sP_k\i P_{s-k+1}\i
+q^{2a}q^{-t^2+(t-1)^2}q^{(t-a)^2-(t-a-1)^2}(q^{2t-2a}-1)P_sP_k\i P_{s-k+1}\i\\&=
q^{2t-2a}P_{2t+1}P_{2a+1}\i P_{2t-2a+1}\i\endalign$$
and (iii) follows.

(iv) {\it If $s\in\ZZ'$, $k\in\ZZ''_{\ge2}$, $\d\in\{1,-1\}$ then 
$N_{s,k}^{*,\d}=2\i P_s(P_k^\d)\i P_{s-k}\i$.}
\nl
This is clear: the orthogonal group of $u$ acts transitively on $S_k^\d\cv$.

We now consider a descending sequence $r_0\ge r_2\ge r_4\ge\do$ of integers such that
$r_0=s$ and $r_{2n}=0$ for large $n$. We denote by $\nu(r_0,r_2,r_4,\do)$ (if $s\in\ZZ'$)  
or by $\nu^\e(r_0,r_2,r_4,\do)$ (if $s\in\ZZ''_{\ge2}$, $\et_u=\e$) the number of 
sequences $U_0\sps U_2\sps U_4\sps\do$ of subspaces of $\cv$ such that for any $n\ge0$,  
$\dim U_{2n}=r_{2n}$ and the quadratic form $u|_{U_{2n}}$ is nondegenerate. We show:

(a) {\it Let $X=\nu(r_0,r_2,r_4,\do)$ if $r_0\in\ZZ'$ and $X=\nu^\e(r_0,r_2,r_4,\do)$ if 
$r_0\in\ZZ''_{\ge2}$ and $\e\in\{1,-1\}$. Then $X$ is a polynomial in $q$ with coefficients
in $\ZZ$ independent of $q$ or $p$.}
\nl
Let $M$ be the number of nonzero terms in the sequence $(r_0,r_2,r_4,\do)$. We have 
$M\ge1$. We argue by induction on $M$. If $M=1$ we have $X=1$. Assume now that $M\ge2$. If 
$r_{2n}=r_{2n+2}>0$ for some $n$ we have $X=\nu(r_0,\do,r_{2n},r_{2n+4},\do)$ if 
$r_0\in\ZZ'$, 
$X=\nu^\e(r_0,\do,r_{2n},r_{2n+4},\do)$ if $r_0\in\ZZ''$; the result follows. Hence we
may assume that the nonzero terms of $r_0,r_2,\do$ are distinct. If $r_{2n}\in\ZZ'$ for 
some $n>0$, we have $X=\nu(r_0,\do,r_{2n},0,0,)\nu(r_{2n},r_{2n+2},\do)$ if $r_0\in\ZZ'$,
$X=\nu^\e(r_0,\do,r_{2n},0,0,)\nu(r_{2n},r_{2n+2},\do)$ if $r_0\in\ZZ''$; in both cases the
induction hypothesis applies to the second factor in the right hand side. Thus we may 
assume that any odd number in the sequence $r_0,r_2,\do$ must appear either as $r_0$ or it
must be followed by $0$. Thus we must consider four cases.

{\it Case 1.} {\it $r_{2n}\in\ZZ''$ for all $n$.} Let $r_{2m}$ be the last nonzero term of
$r_0,r_2,\do$. If $m=0$ we have $X=1$. If $m>0$ we have, using (ii): 
$$\align&X=\sum_{\d_1,\d_2,\do,\d_m\in\{1,-1\}}N_{r_0,r_2}^{\e,\d_1}N_{r_2,r_4}^{\d_1,\d_2}
\do N_{r_{2m-2},r_{2m}}^{\d_{m-1},\d_m}\\&=2^{-m}\sum_{\d_1,\d_2,\do,\d_m\in\{1,-1\}}
P_{r_0}^\e(P_{r_0-r_2}^{\e\d_1})\i(P_{r_2-r_4}^{\d_1\d_2})\i\do
(P_{r_{2m-2}-r_{2m}}^{\d_{m-1}\d_m})\i(P_{r_{2m}}^{\d_m})\i.\endalign$$
This is clearly a polynomial in $q$ with coefficients in $\ZZ[2\i]$. It is also the 
product of a polynomial in $q$ with coefficients in $\ZZ$ with
$$\align&2^{-m}\sum_{\d_1,\d_2,\do,\d_m\in\{1,-1\}}
(q^{(r_0-r_2)/2}-\e\d_1)\i(q^{(r_2-r_4)/2}-\d_1\d_2)\i\do\T\\&
\T(q^{(r_{2m-2}-r_{2m})/2}-\d_{m-1}\d_m)\i(q^{r_{2m}/2}-\d_m)\i\\&=(q^{r_0/2}+\e)
(q^{r_0-r_2}-1)\i(q^{r_2-r_4}-1)\i\do(q^{r_{2m-2}-r_{2m}}-1)\i(q^{r_{2m}}-1)\i\endalign$$
hence is a power series with integer coefficients in $q$. Hence $X\in\ZZ[q]$.

{\it Case 2.} {\it $r_0\in\ZZ''$, $r_{2m}\in\ZZ'$, $r_{2n}\in\ZZ''$ for $0<n<m$ and
$r_{2n}=0$ for $n>m$.} If $m=1$ we have $X=N_{r_0,r_2}^{\e,*}$. If $m>1$ we have 
$$X=\sum_{\d_1,\d_2,\do,\d_{m-1}\in\{1,-1\}}N_{r_0,r_2}^{\e,\d_1}N_{r_2,r_4}^{\d_1,\d_2}\do
N_{r_{2m-4},r_{2m-2}}^{\d_{m-2},\d_{m-1}}N_{r_{2m-2},r_{2m}}^{\d_{m-1},*}.$$
As in case 1 we see, using (ii) and (i) that $X\in\ZZ[q]$.

{\it Case 3.} {\it $r_0\in\ZZ'$, $r_{2n}\in\ZZ''$ for $n>0$.} Let $r_{2m}$ be the last 
nonzero term of $r_0,r_2,\do$. We have 
$$X=\sum_{\d_1,\d_2,\do,\d_m\in\{1,-1\}}N_{r_0,r_2}^{*,\d_1}N_{r_2,r_4}^{\d_1,\d_2}\do 
N_{r_{2m-2},r_{2m}}^{\d_{m-1},\d_m}.$$
As in case 1 we see, using (ii) and (iv) that $X\in\ZZ[q]$.

{\it Case 4.} {\it $r_0\in\ZZ'$ and for some $m>0$, $r_{2m}\in\ZZ'$, $r_{2n}\in\ZZ''$ for 
$0<n<m$ and $r_{2n}=0$ for $n>m$.} If $m=1$ we have $X=N_{r_0,r_2}^{*,*}$, see (iii). If 
$m>1$ we have
$$X=\sum_{\d_1,\d_2,\do,\d_{m-1}\in\{1,-1\}}N_{r_0,r_2}^{*,\d_1}N_{r_2,r_4}^{\d_1,\d_2}\do
N_{r_{2m-4},r_{2m-2}}^{\d_{m-2},\d_{m-1}}N_{r_{2m-2},r_{2m}}^{\d_{m-1},*}.$$
As in case 1 we see, using (i), (ii) and (iv) that $X\in\ZZ[q]$.
This completes the proof of (a).

\mpb

For any $m\in\ZZ''_{\ge0}$ we set

$R_m=q^{m^2/4)}(q^2-1)(q^4-1)\do(q^m-1)$.
\nl
Let $\cv'$ be a $\kk$-vector space with a nondegenerate symplectic form 
$(,):\cv'\T\cv'@>>>\kk$.

We consider a descending sequence $r_1\ge r_3\ge r_5\ge\do$ of even integers such that
$r_1=\dim\cv'$ and $r_{2n+1}=0$ for large $n$. We denote by $\nu'(r_1,r_3,r_5,\do)$ the 
number of sequences $U_1\sps U_3\sps U_5\sps\do$ of subspaces of $\cv'$ such that for any 
$n\ge0$, $\dim U_{2n+1}=r_{2n+1}$ and the symplectic form $(,)|_{U_{2n+1}}$ is 
nondegenerate. Clearly,

(b) {\it $\nu'(r_1,r_3,r_5,\do)=R_{r_1}R\i_{r_1-r_3}R\i_{r_3-r_5}\do$ is a 
polynomial in $q$ with coefficients in $\ZZ$ independent of $q$ or $p$.}

\subhead 1.3\endsubhead
Let $\bV\in\bacc$. Define $(f_a)_{a\in\ZZ}$ by $f_a=\dim\bV^a$. A quadratic form 
$\bQ:\bV@>>>\kk$ with associated symmetric bilinear form $\la,\ra$ is said to be compatible
with the grading if $\bQ|_{\bV^a}=0$ for $a\ne0$ and $\la\bV^a,\bV^{a'}\ra=0$ for
$a+a'\ne0$. In this subsection we fix such a $\bQ$ (compatible with
the grading) and we assume that it is nondegenerate. Let $\bR=\Rad_{\la,\ra}$. Then 
$\bR\sub\bV^0$, $Q|_{\bV^0}:\bV^0@>>>\kk$ is nondegenerate and $\la,\ra$ restricts to a 
perfect pairing $\bV^{-a}\T\bV^a@>>>\kk$ for any $a\ge1$. Hence $f_a=f_{-a}$ for all $a$. 
Let 
$$\align E^2\bV&=\{T\in\Hom(\bV,\bV);T(\bV^a)\sub\bV^{a+2}\text{ for any }a\in\ZZ;
\la Tx,y\ra\\&+\la x,Ty\ra=0\text{ for all }x,y\in\bV;\la x,Tx\ra=0\text{ for any }
x\in\bV^{-1}\}.\endalign$$
If $T\in E^2\bV$ we have 

(a) $\la x,T^ax\ra=0$ for any $a\in ZZ'_{\ge1}$ and any $x\in\bV^{-a}$.
\nl
Indeed, we have $a=2a'+1$ with $a'\ge0$ and $x'=T^{a'}x\in\bV^{-1}$ satisfies
$\la x',Tx'\ra=0$. Thus 
$$0=\la T^{a'}x,T^{a'+1}x\ra=(-1)^{a'}\la x,T^a\ra$$
and (a) follows.

Let $E^2_*\bV$ be the set of all $T\in E^2\bV$ such that 

(i) for any $a\in\ZZ''_{\ge0}$, 
the quadratic form $\bQ_a:\bV^{-a}@>>>\kk$, $x\m\bQ(T^{a/2}x)$ is
nondegenerate;

(ii) for any $a\in ZZ'_{\ge1}$, the symplectic form $\bV^{-a}\T\bV^{-a}@>>>\kk$, 
$x,y\m\la x,T^ay\ra$ is nondegenerate.
\nl
Note that (i) is automatic for $a=0$. 

For $\tN\in E^2\bV$ and $a\ge0$ we set $K_a^T=\ker(T^a:\bV^{-a}@>>>\bV^a)$. If $a\ge1$ then
$K_a^T=\Rad_{\la,\ra_a}$ where $\la,\ra_a:\bV^{-a}\T\bV^{-a}@>>>\kk$ is 
$x,y\m\la x,T^ay\ra$. Hence condition (ii) is equivalent to the condition that $K_a^T=0$ 
for $a\in ZZ'_{\ge1}$. 
If $a\in\ZZ''_{\ge0}$ then the symmetric bilinear form associated to $\bQ_a$
is $(-1)^{a/2}\la,\ra_a$. Hence condition (i) is equivalent to the condition that 
$\bQ_a:K_a^T@>>>\kk$ is injective for $a\in\ZZ''_{\ge2}$. (It implies that for such $a$ we
have $K_a^T=0$ unless $p=2$ and $f_{-a}\in\ZZ'$, when $\dim K_a^T=1$.) 

This discussion shows that if $p\ne2$ an element $T\in E^2\bV$ belongs to $E^2_*\bV$ if and
only if $T^a:\bV^{-a}@>>>\bV^a$ is an isomorphism for any $a\ge0$.

Returning to the general case we reformulate the conditions (i),(ii) for an element 
$T\in E^2\bV$ to be in $E^2_*\bV$ as follows:

(i${}'$) for any $a\in\ZZ''_{\ge0}$, the map $T^{a/2}:\bV^{-a}@>>>\bV^0$ is injective and its
image $I_a^T$ is such that $\bQ|_{I_a^T}$ is a nondegenerate quadratic form;

(ii${}'$) the symplectic form $\o_T(x,y)=\la x,Ty\ra$ on $V^{-1}$ is nondegenerate; for any
$a\in ZZ'_{\ge1}$, 
the map $T^{(a-1)/2}:\bV^{-a}@>>>\bV^{-1}$ is injective and its image $I_a^T$
is such that $\o_T|_{I_a^T}$ is a nondegenerate symplectic form.
\nl
Clearly if (i${}'$) holds then (i) holds. If (i) holds and 
$x\in\ker(T^{a/2}:\bV^{-a}@>>>\bV^0)$ (with $a\in\ZZ''_{\ge2}$) then $T^ax=0$ hence 
$x\in K_a^T$. We have also $\bQ_a(x)=0$. Since $\bQ_a:K_a^T@>>>\kk$ is 
injective we see that $x=0$. We see that (i${}'$) holds.

For $a\in ZZ'_{\ge1}$ and $x,y\in\bV^{-a}$ we have 

$\la x,y\ra_a=(-1)^{(a-1)/2}\la T^{(a-1)/2}x,T^{(a-1)/2}y\ra_1$.
\nl
Hence if (ii${}'$) holds 
then (ii) holds. If (ii) holds and $x\in\ker(T^{(a-1)/2}:\bV^{-a}@>>>\bV^{-1})$ (with 
$a\in ZZ'_{\ge1}$) then $T^ax=0$ hence $x\in K_a^T$ and $x=0$. We see that (ii${}'$) holds.

If $T\in E^2_*\bV$ then clearly,

$\bV^0=I_0^T\sps I_2^T\sps I_4^T\sps\do$ and $\bV^{-1}=I_1^T\sps I_3^T\sps I_5^T\sps\do$.
\nl
We see that if $E^2_*\bV\ne\em$ then

(b) $f_0\ge f_{-2}\ge f_{-4}\ge\do$ and $f_{-1}\ge f_{-3}\ge f_{-5}\ge\do$,

(c) $f_a\in\ZZ''$ if $a\in\ZZ'$.
\nl
We say that $(f_a)_{a\in\ZZ}$ is {\it admissible} if it satisfies (b),(c). In the remainder
of this subsection we assume that $(f_a)_{a\in\ZZ}$ is admissible. Let $R^0$ be the set of
all sequences $U_0\sps U_2\sps U_4\sps\do$ of subspaces of $\bV^0$ such that 
$\dim U_a=f_{-a}$ and $\bQ|_{U_a}$ is a nondegenerate quadratic form for $a=0,2,4,\do$. Let
$R^{-1}$ be the set of all pairs $(\o,(U_1,U_3,U_5,\do))$ where $\o$ is a nondegenerate 
symplectic form on $\bV^{-1}$ and $U_1\sps U_3\sps U_5\sps\do$ are subspaces of $\bV^{-1}$
such that $\dim U_a=f_{-a}$ and $\o|_{U_a}$ is a nondegenerate symplectic form for 
$a=1,3,5,\do$. Clearly $R^0\ne\em,R^{-1}\ne\em$. Define $\ps:E^2_*\bV@>>>R^0\T R^{-1}$ by 

$T\m((I_0^T,I_2^T,I_4^T,\do),(\o_T,(I_1^T,I_3^T,I_5^T,\do)))$. 
\nl
Clearly, 

(d) {\it the fibre of $\ps$ at $((U_0,U_2,U_4,\do),(\o,(U_1,U_3,U_5,\do)))\in R^0\T R^{-1}$
can be identified with $\prod_{a\ge2}\Iso(\bV^{-a},U_a)$};
\nl
here $\Iso(\bV^{-a},U_a)$ is the set of vector space isomorphisms $\bV^{-a}@>>>U_a$.

We now assume that $\kk,q$ are as in 1.2. For any $m\ge0$ we set

$A_m=q^{1/2(n^2-n)}(q-1)(q^2-1)\do(q^m-1)$.
\nl
If $f_0\in\ZZ''_{\ge2}$ let $\e=\et_{\bQ|_{\bV^0}}$. From (d) we see that 

(e) $|E^2_*\bV|=\prod_{a\ge1}A_{f_a}\x\nu'(f_1,f_3,f_5,\do)(R_{f_1})\i$
\nl
where $\x=\nu(f_0,f_2,f_4,\do)$ if $f_0\in\ZZ'$; $\x=\nu^\e(f_0,f_2,f_4,\do)$ if 
$f_0\in\ZZ''_{\ge2}$; $\x=1$ if $f_0=0$.

\subhead 1.4\endsubhead
Let $X^*=(X^{\ge a})_{a\in\ZZ}$ be a sequence of subspaces of $V$ such that 
$X^{\ge a+1}\sub X^{\ge a}$ for any $a$, $X^{\ge a}=0$ for some $a$, $X^{\ge a}=V$ for some
$a$ that is a {\it filtration} (see \cite{\LU, 2.2}) of $V$. We say that $X^*$ is a
$Q$-filtration of $V$ if for any $a\ge1$ we have

(a) $Q|_{X^{\ge a}}=0$ and $X^{\ge1-a}=(X^{\ge a})^\pe$.
\nl
Then for any $a\le0$ we have:

(b) {\it $(X^{\ge a})^\pe=X^{\ge1-a}+R$ (direct sum) and
$X^{\ge1-a}=(X^{\ge a})^\pe\cap Q\i(0)$.}
\nl
The first equality follows by applying ${}^\pe$ to both sides of 
$X^{\ge a}=(X^{\ge1-a})^\pe$, see (a). If $x\in V^{\ge1-a}\cap R$ then by (a) we have 
$Q(x)=0$. Since $Q:R@>>>\kk$ is injective, we have $x=0$. Thus, $V^{\ge1-a}\cap R=0$. The
second equality in (b) follows from the first equality and (a). This proves (b).

The proof of the following result is standard.

(c) {\it We can find a direct sum decomposition $V=\op_{a\in\ZZ}X^a$ such that 
$X^{\ge a}=\op_{a';a'\ge a}X^a$ for all $a$ and $\la X^a,X^{a'}\ra=0$ for $a+a'\ne0$,
$Q|_{X^a}=0$ for $a\ne0$.}
\nl
For any $Q$-filtration $X^*$ of $V$ we set 
$$E^{\ge2}X^*=\{N\in\tcm_Q;NX^{\ge a}\sub X^{\ge a+2}\text{ for all }a\in\ZZ\}.$$
We show:

(d) {\it If $x\in X^{\ge-1}$ then $\la x,Nx\ra=0$.}
\nl
It is enough to show that $Q(Nx)=0$. Since $Nx\in X^{\ge1}$, this follows from (a).

\subhead 1.5\endsubhead
Let $X^*$ be a $Q$-filtration of $V$. We set $\gr^aX^*=X^{\ge a}/X^{\ge a+1}$. Then 
$\gr X^*=\op_{a\in\ZZ}\gr^aX^*\in\bacc$. Let $f_a=\dim\gr^aX^*$. Define a quadratic form 
$\bQ:\gr X^*@>>>\kk$ by $\bQ(x)=Q(\dx_0)+\sum_{a\ge1}\la\dx_{-a},\dx_a\ra$ where 
$x=\sum_ax_a$ (with $x_a\in\gr^aX^*$) and $\dx_b$ is a representative of $x_b$ in 
$X^{\ge b}$. The symmetric bilinear form associated to $\bQ$ is 
$\sum_ax_a,\sum_bx'_b\m\sum_{a+b=0}\la\dx_a,\dx'_b\ra$ where $x_a\in\gr^aX^*$,
$x'_b\in\gr^bX^*$ and $\dx_a,\dx'_b$ are representatives of $x_a,x'_b$ in $X^{\ge a}$,
$X^{\ge b}$. This form is denoted again
by $\la,\ra$; its radical is the image of $R$ under $X^{\ge0}@>>>\gr^0X^*$. It follows that
$\bQ$ is nondegenerate. It is clearly compatible with the grading. Hence $f_a=f_{-a}$ for
any $a$.

Now let $N\in E^{\ge2}X^*$. For any $a$, $N$ restricts to a linear map
$X^{\ge a}@>>>X^{\ge a+2}$ and $X^{\ge a+1}@>>>X^{\ge a+3}$ hence it induces a linear map
$\gr^aX^*@>>>\gr^{a+2}X^*$. Taking the direct sum over $a$ of these linear maps we obtain a
linear map $\bN:\gr X^*@>>>\gr X^*$. We have $\bN\in E^2\gr X^*$.
(If $x\in\gr^aX^*,y\in\gr^bX^*$ and $\dx,\dy$ are representatives for $x,y$ in 
$X^{\ge a},X^{\ge b}$, the sum $\la\bN x,y\ra+\la x+\bN y\ra$ is $0$ unless $a+b+2=0$ in 
which case it $\la N\dx,\dy\ra+\la\dx+N\dy\ra=-\la N\dx,N\dy\ra$ which is again zero since
$N\dx\in X^{\ge a+2}$, $N\dy\in X^{\ge b+2}=X^{\ge-a}$ and
$\la X^{\ge a+2},X^{\ge-a}\ra=0$. If $x,dx$ are as above and $a=-1$, we have
$\la x,\bN x\ra=\la\dx,N\dx\ra=-Q(N\dx)$ and this is $0$ since $N\dx\in X^{\ge1}$ and
$Q|_{X^{\ge1}}=0$.) Thus we have a well defined map
$$\Ph:E^{\ge2}X^*@>>>E^2\gr X^*,\qua N\m\bN.$$
Let 
$$d=\sum_{a<a';-a-a'\ge3}f_af_{a'}+\sum_{a;-2a\ge4}f_a(f_a-1)/2.$$
We show:

(a) {\it $\Ph$ is an iterated affine space bundle with fibres of dimension $d$.}
\nl
Let $V=\op_aX^a$ be as in 1.4(c). For any integer $k\ge2$ let $Z_k$ be the set of all 
collections $(N^a_b)_{a,b\in\ZZ;2\le b-a\le k}$ where $N^a_b:X^a@>>>X^b$ are linear maps 
that satisfy

(i) $\la x_a,N^{a'}_{-a}x'_{a'}\ra+\la N^a_{-a'}x_a,x'_{a'}\ra=
-\sum_{b;-a'-2\ge b\ge a+2}\la N^a_bx_a,N^{a'}_{-b}x'_{a'}\ra$
\nl
for any $a\ne a'$ such that $2\le-a-a'\le k$ and any $x_a\in X^a,x'_{a'}\in X^{a'}$;

(ii) $\la x_a,N^a_{-a}x_a\ra
=-Q(N^a_0x_a)-\sum_{b<0;-a-2\ge b\ge a+2}\la N^a_bx_a,N^a_{-b}x_a\ra$
\nl
for any $a$ such that $4\le-2a\le k0$ and any $x_a\in X^a$;

(iii) $\la x_{-1},N^{-1}_1x_{-1}\ra=0$
\nl
for any $x_{-1}\in X^{-1}$.

For large $k$ we may identify $Z_k=E^{\ge2}X^*$ by $(N^a_b)\m N$,
$Nx_a=\sum_{b;b\ge a+2}N^a_bx_a$ with $x_a\in X^a$. Moreover we may identify
$Z_2=E^2\gr X^*$ in an obvious way. We have obvious maps $Z_2@<<<Z_3@<<<Z_4@<<<\do$. These
maps eventually become the identity map of $E^{\ge2}X^*$; their composition may be 
identified with $\Ph$. It is enough to show that for any $k\ge3$ the obvious map 
$Z_k@>>>Z_{k-1}$ is an affine space bundle with fibres of dimension 
$$d_k=\sum_{a<a';-a-a'=k}f_af_{a'}+\sum_{a;-2a=k}f_a(f_a-1)/2.$$
We shall prove only that the fibre of this map at any given point of $Z_{k-1}$ is an affine
space of dimension $d_k$. This fibre may be identified with the set of all collections 
$(N^a_b)_{a,b\in\ZZ;b-a=k}$ where $N^a_b:X^a@>>>X^b$ are linear maps that satisfy (i) with
$-a-a'=k$, (ii) with $-2a=k$. (In these equations the right hand sides 
involve only the coordinates of the given point in $Z_{k-1}$.) In the equation (i) with 
$a\ne a',-a-a'=k$ each of $N^{a'}_{-a},N^a_{-a'}$ determines the other. So the solutions of
this equation form an affine space of dimension $f_af_{a'}$. If $k\in\ZZ'$ there are no 
further equations. If $k\in\ZZ''$, in the equation (ii) with $-2a=k$ the right hand side is
a known quadratic form on $X^a$ and the solutions $N^a_{-a}$ form an affine space of 
dimension $f_a(f_a-1)/2$. This completes the proof of (a). 

We show:

(b) {\it $1+E^{\ge2}X^*\sub SO_Q$.}
\nl
To prove this we may assume that $|\kk|=\iy$. Since $1\in 1+E^{\ge2}X^*\sub O_Q$ it is 
enough to show that $E^{\ge2}X^*$ is irreducible. Using (a) we see that it is enough to 
show that $E^2\gr X^*$ is irreducible. From the definitions we see that $E^2\gr X^*$ is an
affine space of dimension $\sum_{a<a';-a-a'=2}f_af_{a'}+f_{-1}(f_{-1}-1)/2$. This proves 
(b).

\mpb

We set $E^{\ge2}_*X^*=\Ph\i(E^2_*\gr X^*)$. From (a) we see that:

(c) {\it $E^{\ge2}_*X^*$ is (via $\Ph$) an iterated affine space bundle over $E^2_*\gr X^*$
with fibres of dimension $d$.}
\nl
Using this and results in 1.3 we see that if $E^{\ge2}_*X^*\ne\em$ then 
$(f_a)_{a\in\ZZ}$ is admissible.

In the remainder of this subsection we assume that $\kk,q$ are as in 1.2 and that 
$(f_a)_{a\in\ZZ}$ is admissible. From (c) we see that

(d) $|E^{\ge2}_*X^*|=q^d|E^2_*\gr X^*|$.
\nl
We denote $|E^{\ge2}_*X^*|$ by $\fB_{(f_a)}$ if $D\in\ZZ'$; by $\fB^\e_{(f_a)}$ if 
$D\in\ZZ''_{\ge2}$ and $\et_Q=\e$. From (d), 1.3(e), 1.2(a) and 1.2(b) we see that:

(e) {\it Let $Y=\fB_{(f_a)}$ if $f_0\in\ZZ'$; $Y=\fB^\e_{(f_a)}$ if $f_0\in\ZZ''$, $\sum_af_a>0$,
$\e\in\{1,-1\}$. Then $Y$ is a polynomial in $q$ with coefficients in $\ZZ$ independent of 
$q$ or $p$.}
\nl
Note that if $D\in\ZZ''_{\ge2}$ we have $\et_Q=\et_{\bQ}$.

\subhead 1.6\endsubhead
In the remainder of this section we assume that $D\ge2$. Let $\tcf^D$ be the set of all 
collections $(f_a)_{a\in\ZZ}$ of natural numbers such that $\sum_af_a=D$, $f_a=f_{-a}$ for
all $a$ and the admissibility conditions 1.3(b),(c) hold. When $D\in\ZZ''$, let $J$ be the
set of $SO(V)$-orbits on $\cj_Q$ (if $|\kk|=\iy$) or on $\cj_{\ti Q}$ for $\ti Q$ as in 
1.1 (if $|\kk|<\iy$); note that $|J|=2$; in this case let 

$\tcf_0^D=\{(f_a)\in\tcf^D;f_0=0\}$, $\tcf_1^D=\{(f_a)\in\tcf^D;f_0>0\}$,

$\cf^D=\tcf_1^D\sqc(J\T\tcf_0^D)$.
\nl
When $D\in\ZZ'$ we set $\cf^D=\tcf^D$. For 
$\ph\in\cf^D$ of the form $(f_a)$ or $(j,(f_a))$ (where $j\in J$), let $\bar\cy_\ph$ be 
the set of all $Q$-filtrations  $X^*$ of $V$ such that $\dim\gr^aX^*=f_a$ for all $a$ and 
(in the case where $\ph=(j,(f_a))$), $X^0=X^1\in j$; let $\cy_\ph$ be the set of all pairs
$(X^*,N)$ such that $X^*\in\bar\cy_\ph$ and $N\in E^{\ge2}_*X^*$. Note that $\bar\cy_\ph$ 
is a partial flag manifold of $SO_Q$. Moreover, the obvious map $\cy_\ph@>>>\bar\cy_\ph$ 
has fibres $E^{\ge2}_*X^*$ which are smooth (if $|\kk|=\iy$) since $E^{\ge2}_*X^*$ is open
in the affine space $E^{\ge2}X^*$. We see that $\cy_\ph$ is naturally a smooth variety (if
$|\kk|=\iy$).

Define $\Ps:\sqc_{\ph\in\cf^D}\cy_\ph@>>>\{g\in SO_Q;g\text{ unipotent}\}$ by 
$(X^*,N)\m1+N$ for $(X^*,N)\in\cy_\ph$. This is well defined by 1.5(b).

\proclaim{Theorem 1.7}In the setup of 1.6, $\Ps$ is a bijection.
\endproclaim
An equivalent statement is:

(a) {\it For any unipotent element $g\in SO_Q$ there is a unique $Q$-filtration $X^*$ of
$V$ such that $g-1\in E^{\ge2}_*X^*$.}
\nl
It is enough to prove this in the case where $|\kk|=\iy$. The proof is given in 2.12, 3.3.

\mpb

For any $\ph\in\cf^D$ let $\Xi_\ph=\Ps(\cy_\ph)$. The theorem shows that the sets 
$\Xi_\ph$ form a partition of the variety of unipotent elements in $SO_Q$ and that for any
$\ph\in\cf^D$, $\Ph$ restricts to a bijection $\cy_\ph@>>>\Xi_\ph$. One can show that if
$|\kk|=\iy$ this is an isomorphism of the smooth variety $\cy_\ph$ with a subvariety of the
variety of unipotent elements of $SO_Q$. Note that if $|\kk|=\iy$ then each $\cy_\ph$ is 
nonempty hence each $\Xi_\ph$ is nonempty. (If $|\kk|<\iy$ then $\Xi_\ph$ is nonempty 
unless $D\in\ZZ''$, $\et_Q=-1$, $\ph\in J\T\tcf_0^D$ in which case it is empty.) 
In the case where $|\kk|=\iy$ and $p\ne2$, the subsets $\Xi_\ph$ are exactly the unipotent
conjugacy classes in $SO_Q$. If $p=2$, the subsets $\Xi_\ph$ are unions of unipotent 
conjugacy classes in $SO_Q$ but not necessarily single conjugacy classes; see Section 4.

\subhead 1.8\endsubhead
In the setup of 1.6 we assume that $\kk,q$ are as in 1.2. We assume that Theorem 1.7 holds.
Let $\ph\in\cf^D$. Assume that if $D\in\ZZ''$, $\et_Q=-1$, then $\ph\in\tcf_1^D$. We have

(a) $|\cy_\ph|=|\bar\cy_\ph||E^{\ge2}_*X^*|$
\nl
for any $X^*\in\bar\cy_\ph$. We denote $|\cy_\ph|$ by $\ti\fB_\ph$ if $D\in\ZZ'$; by 
$\ti\fB^\e_\ph$ if $D\in\ZZ''$ and $\et_Q=\e$. From (a) and 1.5(e) we see that:

(b) {\it Let $\tY=\ti\fB_\ph$ if $D\in\ZZ'$; $\tY=\ti\fB^\e_\ph$ if $D\in\ZZ''$ and 
$\et_Q=\e$. Then $\tY$ is a polynomial in $q$ with coefficients in $\ZZ$ independent of $q$
or $p$.}

\subhead 1.9\endsubhead
Let $N\in\tcm_Q-\{0\},e=e_N$. We have $e\ge2$. Let E be a complement to $\ker N^{e-1}$ in 
$V$. We show

(a) {\it The $(-1)^{e-1}$-symmetric bilinear form $(,):E\T E@>>>\kk$ given by 
$(x,y)=\la x,N^{e-1}y\ra$ is nondegenerate.}
\nl
Assume that $v_0\in E,\la v_0,N^{e-1}E\ra=0$. Then $\la N^{e-1}v_0,E\ra=0$. Also, 

$\la N^{e-1}v_0,\ker N^{e-1}\ra=\pm\la v_0,N^{e-1}\ker N^{e-1}\ra=0$.
\nl
Thus $\la N^{e-1}v_0,E+\ker N^{e-1}\ra=0$ that is $\la N^{e-1}v_0,V\ra=0$ hence 
$N^{e-1}v_0\in R$. Using 1.1(a) we deduce $N^{e-1}v_0=0$. We see that
$v_0\in\Ker N^{e-1}\cap E=0$. This proves (a).

We show:

(b) {\it the linear map $E^{\op e}@>>>V$ given by
$(v_0,v_1,\do,v_{e-1})\m v_0+Nv_1+\do+N^{e-1}v_{e-1}$ is injective and $\la,\ra$ is 
nondegenerate on its image $W:=E+NE+\do+N^{e-1}E$.}
\nl
Assume that $v=v_0+Nv_1+\do+N^{e-1}v_{e-1}$ with $v_i\in E$ satisfies
$0=\la v,E\ra=\la v,NE\ra=\do=\la v,N^{e-1}E\ra$. We must show that $v_k=0$ for 
$k\in[0,e-1]$. We argue by induction on $k$. Assume that $k\in[0,e-1]$ and that $v_{k'}=0$
for $k'<k$. From $0=\la v,N^{e-1-k}E\ra$ we get $\la N^kv_k,N^{e-1-k}E\ra=0$ hence 
$(v_k,E)=0$ hence $v_k=0$ (using (a)). This proves (b).

We show:

(c) {\it we have $V=W\op Y$ ($W$ as in (b), $Y=W^\pe$), $R\sub Y$ and $Y$ is $N$-stable; 
moreover $N^{e-1}Y=0$.}
\nl
The first assertion follows from (b). Now let $x\in Y$. We can write $x=x'+x''$ where
$x'\in E$, $x''\in\ker N^{e-1}$. We have $N^{e-1}x=N^{e-1}x'\in W\cap Y$ hence 
$N^{e-1}x=0$. This proves (c).

\head 2. The case $p=2$\endhead
\subhead 2.1\endsubhead
In this section we assume that $p=2$. In this case the symmetric bilinear form $\la,\ra$ 
associated to $Q$ is symplectic. Let $N\in\tcm_Q$.

For $i\ge1$ we define $\l_i^N:\ker N^i@>>>\kk$ by $x\m\la x,N^{i-1}x\ra^{1/2}$. Note that 
$\l_i^N=0$ if $i\in\ZZ'$ and $\l_i^N$ is linear for any $i$. 

For $i\ge1$ we define $\e_i^N\in\{0,1\}$ by $\e_i^N=0$ if $\l_i^N=0$, $\e_i^N=1$ if 
$\l_i^N\ne0$. Note that $\e_i^N=0$ if $i\in\ZZ'$. 

\subhead 2.2\endsubhead
Let $N\in\tcm_Q-\{0\},e=e_N$. We have $e\ge2$. Let $\l=\l_e^N$. Define a subspace $L=L_N$ 
of $V$ by 

$L=N^{e-1}V$ if $\l=0$, 

$L=(\ker\l)^\pe$ if $\l\ne0,R=0$,

$L=\{x\in(\ker\l)^\pe;q(x)=0\}$ if $\l\ne0,R\ne0$.
\nl
If $\l\ne0$ then $\ker\l$ is a hyperplane in $V$ containing $R$; thus $L$ is a line if
$R=0$. If $\l\ne0,R\ne0$ then $K:=(\ker\l)^\pe$ is a two dimensional subspace containing 
$R$ so that $\la,\ra$ is zero on $K$. Hence $q:K@>>>\kk$ is a group homomorphism which
restricts to a group isomorphism $R@>>>\kk$ hence induces a group isomorphism $K/L@>>>\kk$.
We see that $K=L\op R$. Thus $L$ is a line complementary to $R$ in $K$ and $\ker\l=L^\pe$.

We show:

(a) $L\sub L^\pe$.
\nl
Assume first that $\l=0$. By 1.1(b) it is enough to show that $N^{e-1}V\sub\ker N^{e-1}$
which is clear since $2e-2\ge e$. If $\l\ne0$ then (a) is clear since $L$ is a line.

We show:

(b) {\it if $\l\ne0$, then $L\sub N^{e-1}(L^\pe)$ if and only if $N^{e-1}V\in\cc''$.}
\nl
Define a symmetric bilinear form $(,)$ on $V$ by $(y,y')=\la y,N^{e-1}y'\ra$. We have 
$\Rad_{(,)}=\{y'\in V;N^{e-1}y'\in R\}=\ker N^{e-1}$ (since $N^{e-1}V\cap R=0$). We have 
$\l(y)=(y,y)^{1/2}$. Since $\l|_{\Rad_{(,)}}=0$, we can find $v\in V$ such that 
$(v,y)=\l(y)$ that is $\la N^{e-1}v,y\ra=\l(y)$ for any $y\in V$. Since $\l\ne0$ we have 
$N^{e-1}v\n0$. Clearly, if $y\in\ker\l$ then $\la N^{e-1}v,y\ra=0$. Thus 
$N^{e-1}v\in(\ker\l)^\pe$.

The symmetric bilinear form $(,)_0$ on $V$ given by $(y,y')_0=(y,y')+\l(y)\l(y')$ is 
symplectic: for $y\in V$ we have $(y,y)_0=(y,y)+\l(y)^2=0$. Since $\l|_{\ker N^{e-1}}=0$ we
see that $R_0:=\Rad_{(,)_0}$ contains $\ker N^{e-1}$. Since $(,)_0$ is symplectic we have 
$\dim R_0=D\mod2$.

Assume now that $N^{e-1}V\in\cc''$. Then $\dim\ker N^{e-1}=D\mod2$. Hence 
$\dim\ker N^{e-1}=\dim R_0\mod2$ that is $R_0/\ker N^{e-1}\in\cc''$. If 
$R_0\ne\ker N^{e-1}$ then $\dim(R_0/\ker N^{e-1})\ge2$. Hence the kernel of the linear map
$R_0/\ker N^{e-1}@>>>\kk$ induced by $\l:R_0@>>>\kk$ is non-zero. Thus there exists 
$u\in R_0$ such that $u\n\ker N^{e-1}$ and $\l(u)=0$. For any $y\in V$ we have 
$0=(y,u)_0=(y,u)+\l(u)\l(y)=(y,u)$ so that $u\in\Rad_{(,)}=\ker N^{e-1}$. This 
contradiction shows that $R_0=\ker N^{e-1}$. Since $N^{e-1}v\ne0$ we have 
$v\n\in\ker N^{e-1}$ hence $v\n R_0$. Thus there exists $v'\in V$ such that $(v,v')_0\ne0$
that is, $(v,v')+(v,v)(v,v')\ne0$. It follows that $1+(v,v)\ne0$. Moreover we have 
$0=(v,v)_0=(v,v)+(v,v)(v,v)$ so that $(v,v)(1+(v,v))=0$. Since $1+(v,v)\ne0$ we deduce that
$(v,v)=0$ that is $\la v,N^{e-1}v\ra=0$. We have 
$Q(N^{e-1}v)=\la N^{e-2}v,N^{e-1}v\ra=\la v,N^{2e-3}v\ra$. If $e\ge3$ this is $0$ since 
$2e-3\ge e$; if $e=2$ this is equals $\la v,N^{e-1}v\ra$ which is again $0$. We see that
$N^{e-1}v$ is a non-zero vector in $\ker(Q:(\ker\l)^\pe@>>>\kk)$. Hence 
$N^{e-1}v\in L-0$. Since $L$ is a line we see that $L$ is spanned by $N^{e-1}v$. To show
that $L\sub N^{e-1}(L^\pe)$ it is enough to show that $v\in L^\pe$ or that $\la v,L\ra=0$ 
or that $\la v,N^{e-1}v\ra=0$. But this equality is already known. We see that
$L\sub N^{e-1}(L^\pe)$.

Conversely, assume that $L\sub N^{e-1}(L^\pe)$. Let $x\in L-0$. We have $x=N^{e-1}t$ with
$\la t,L\ra=0$. In particular, $\la t,x\ra=0$ that is $\l(t)=0$. Since $x\in L$ we have 
$z\in\ker\l\imp\la x,z\ra=0$. Hence there exists $c\in\kk$ such that $\la x,z\ra=c\l(z)$ 
for all $z\in V$. Since $x\n R$ we have $c\ne0$. Replacing $x$ by a non-zero scalar 
multiple we may assume that $\la x,z\ra=\l(z)$ that is $(t,z)=\l(z)$ for all $z\in V$. Let
$R_1=\Rad_{(,)|_{\ker\l}}$. If $u\in R_1$ then we have $w\in\ker\l\imp(u,w)=0$. Hence there
exists $c'\in\kk$ such that $(u,w)=c'\l(w)$ that is $(u,w)=c'(t,w)$ for any $w\in V$. Hence
$u-c't\in\Rad_{(,)}=\ker N^{e-1}$. We see that $R_1\sub\kk t+\ker N^{e-1}$. The reverse 
inclusion is obvious. Thus, $R_1=\kk t+\ker N^{e-1}$. We see that $(,)$ induces a 
nondegenerate symmetric bilinear form on $\ker\l/(\kk t+\ker N^{e-1})$. This induced form 
is symplectic: the equality $(y,y')_0=(y,y')+\l(y)\l(y')$ with $(,)_0$ symplectic shows 
that $(,)|_{\ker\l}=(,)_0|_{\ker\l}$. We deduce that $\ker\l/(\kk t+\ker N^{e-1})\in\cc''$.
Since $x\ne0$ we have $t\n\ker N^{e-1}$ hence $\dim\ker\l-1-\dim\ker N^{e-1}\in\ZZ''$ that 
is $D-2-(D-\dim N^{e-1}V)\in\ZZ''$. We see that $N^{e-1}V\in\cc''$. This proves (b).

We show:

(c) We have $L\sub N^{e-1}V+R$. More precisely we have $L\sub N^{e-1}V$ except possibly in
the case where $V\in\cc',\l\ne0,e=2,c_e^N\in\ZZ'$.
\nl
If $\l=0$ this is obvious. Assume now that $\l\ne0$. Since $\ker N^{e-1}\sub\ker\l$ we have
$(\ker\l)^\pe\sub(\ker N^{e-1})^\pe$ and using 1.1(b), $(\ker\l)^\pe\sub N^{e-1}V+R$. Since
$L\sub(\ker\l)^\pe$ we have $L\sub N^{e-1}V+R$. If $R=0$ then clearly $L\sub N^{e-1}V$. If
$e\ge3,R\ne0$ we have $Q|_{N^{e-1}V}=0$; indeed for $v\in V$ we have
$Q(N^{e-1}v)=\la N^{e-2}v,N^{e-1}v\ra=\la v,N^{2e-3}v\ra=0$. Hence if $x\in L$ is written
as $x=N^{e-1}v+v'$ with $v\in V,v'\in R$ then $0=Q(x)=Q(v')$. But $Q(v')=0$, $v'\in R$
implies $v'=0$ hence $x=N^{e-1}v$ and $L\sub N^{e-1}V$. If $c_e^N\in\ZZ''$ then we have
$L\sub N^{e-1}V$ by (b). This proves (c).

We show:

(d) $NL\sub L,\qua N(L^\pe)\sub L^\pe$.
\nl
From (c) we see that $NL=0$ (since $N^e=0$ and $NR=0$). Thus the first inclusion holds. The
second inclusion follows from the first.

\mpb

We define a subset $\cm_Q$ of $\tcm_Q$ as follows:

$\cm_Q=\tcm_Q$ if $V\in\cc'$; 

$\cm_Q=\{\tN\in\tcm_Q;\ker\tN\in\cc''\}$ if $V\in\cc''$.
\nl
We show:

(e) {\it If $N\in\cm_Q-\{0\}$, then $Q|_L=0$.}
\nl
If $R\ne0$ the result is obvious. Now assume that $R=0$. If $e\ge3$ we have 
$Q|_{N^{e-1}V}=0$ as in the proof (c). The same argument shows that $Q|_{N^{e-1}V}=0$ if 
$\l=0$. If $\l=0$ or $e\ge3$ we have $L\sub N^{e-1}V$, see (c), hence $Q|_L=0$. Thus we may
assume that $\l\ne0,e=2$. Since $\ker N\in\cc''$ and $V\in\cc''$ we have $NV\in\cc''$ that
is $N^{e-1}V\in\cc''$. As in the proof of (b) we see that $L$ is spanned by $N^{e-1}v$ ($v$
as in that proof) which is contained in $\ker(Q:(\ker\l)^\pe@>>>\kk)$; the result follows.

\subhead 2.3\endsubhead
Let $N\in\cm_Q-\{0\},e=e_N,\l=\l_e^N,L=L_N$. We have $e\ge2$. We set $V'=L^\pe/L$, see 
2.2(a). From 2.2(d) we see that $N$ induces a (nilpotent) endomorphism of $V'$, denoted by
$N'$. Let $e'=e_{N'}$. If $\l=0$ we have $e'\le e-1$; if $\l\ne0$ we have $e'\le e$. Define
a quadratic form $Q':V'@>>>\kk$ by $Q'(x')=Q(x)$ where $x$ is a representative of 
$x'\in V'$ in $L^\pe$. (To see that $Q'$ is well defined we use 2.2(e).) The symplectic 
form associated to $Q'$ is $\la x',y'\ra'=\la x,y\ra$ where $x',y'\in V'$ and $x,y$ are 
representatives of $x',y'$ in $L^\pe$. Its radical is

$\{x\in L^\pe;\la x,L^\pe\ra=0\}/L=(L+R)/L$. 
\nl
It follows that $Q'$ is nondegenerate. We have $N'\in\tcm_{Q'}$. We show:

(a) $N'\in\cm_{Q'}$.
\nl
To do this we may assume that $V\in\cc''$. Then $R=0$. It is enough to show that 
$\dim\ker N'=\dim\ker N\mod2$. We have $\ker N'=\{x\in L^\pe;Nx\in L\}/L$. We have an exact
sequence 
$$0@>>>L@>a>>\ker N@>b>>\{x\in L^\pe;Nx\in L\}/L@>c>>L\cap N(L^\pe)@>>>0$$
where $a$ is the inclusion, $b$ is induced by the inclusion $\ker N\sub L^\pe$ and $c$ is
induced by $x\m Nx$. Thus $\dim L-\dim\ker N+\dim\ker N'-\dim(L\cap N(L^\pe))=0$ and it is
enough to show that $\dim L=\dim(L\cap N(L^\pe))\mod2$. 

Assume first that $\l\ne0$. We show that $L\sub N(L^\pe)$ (this implies that
$\dim L=\dim(L\cap N(L^\pe))$). If $e=2$ this follows from 2.2(b). (We have $NV\in\cc''$ 
since $\ker N\in\cc''$ and $V\in\cc''$.) Assume now that $e\ge3$. Since $R=0$ we see from 
2.2(c) that $L$ is spanned by $N^{e-1}x$ for some $x\in V$. It is enough to show that 
$N^{e-2}x\in L^\pe$; since $e\ge3$ it is enough to show that $NV\sub L^\pe$; this is clear
since $L^\pe$ is $N$-stable of codimension $1$. 

Assume next that $\l=0$. We must show that 

$\dim N^{e-1}V=\dim(N^{e-1}V\cap N(\ker N^{e-1}))\mod2$.
\nl
If $e\ge3$ we have $N^{e-1}V\sub N(\ker N^{e-1})$ hence 
$\dim N^{e-1}V=\dim(N^{e-1}V\cap N(\ker N^{e-1}))$. If $e=2$ we have $N(\ker N^{e-1})=0$ 
and it is enough to show that $NV\in\cc''$. But this follows from $\ker N\in\cc''$ and 
$V\in\cc''$. This completes the proof of (a).

\subhead 2.4\endsubhead
Let $N\in\cm_Q-\{0\}$. Let $e,\l,L,V',N',Q'$ be as in 2.2. Properties (i)-(iii) below 
describe the invariants $c_i^{N'},\e_i^{N'}$ of $N'$.

(i) {\it Assume that $\l=0$. We have $c_e^{N'}=0$; $c_{e-2}^{N'}=c_{e-2}^N+c_e^N$ (if 
$e>2$); $c_i^{N'}=c_i^N$ if $i\ne e,i\ne e-2$; $\e_i^{N'}=\e_i^N$ for any $i$.}

(ii) {\it Assume that $\l\ne0$ and $c_e^N\in\ZZ''$. Then $c_e^{N'}=c_e^N-2$; 
$c_{e-1}^{N'}=c_{e-1}^N+2$; $c_i^{N'}=c_i^N$ for $i\ne e,i\ne e-1$; $\e_e^{N'}=0$; 
$\e_i^{N'}=\e_i^N$ for any $i\ne e$.}

(iii) {\it Assume that $\l\ne0$, $c_e^N\in\ZZ'$. Then $c_e^{N'}=c_e^N-1$; 
$c_{e-2}^{N'}=c_{e-2}^N+1$ (if $e>2$); $c_i^{N'}=c_i^N$ for $i\ne e,i\ne e-2$; 
$\e_e^{N'}=0$; $\e_{e-2}^{N'}=1$ (if $e>2$); $\e_i^{N'}=\e_i^N$ for any $i\n\{e,e-2\}$.}
\nl
We have $N'{}^{e-1}V'=(N^{e-1}(L^\pe)+L)/L=N^{e-1}(L^\pe)/(N^{e-1}(L^\pe)\cap L)$. 

If $\l=0$ we have $L^\pe=\ker N^{e-1}$ and $N^{e-1}(L^\pe)=0$. Hence $N'{}^{e-1}V'=0$ and
$c_i^{N'}=0$ for $i\ge e$.

If $\l\ne0$ and $N^{e-1}V\in\cc''$ (that is $c_e^N\in\ZZ''$) we have $L\sub N^{e-1}(L^\pe)$
(see 2.2(b)) so that $N'{}^{e-1}V'=N^{e-1}(L^\pe)/L$. If $\l\ne0$ and $N^{e-1}V\in\cc'$ 
(that is $c_e^N\in\ZZ'$) we have $L\not\sub N^{e-1}(L^\pe)$ (see 2.2(b)) and, since $L$ is
a line, we have $N^{e-1}(L^\pe)\cap L=0$ and $N'{}^{e-1}V'=N^{e-1}(L^\pe)$. We have an 
exact sequence 

$0@>>>\ker N^{e-1}\cap L^\pe@>>>L^\pe@>N^{e-1}>>N^{e-1}(L^\pe)@>>>0$.
\nl
Now $\ker N^{e-1}\sub\ker\l=L^\pe$ hence this exact sequence becomes

$0@>>>\ker N^{e-1}@>>>\ker\l@>>>N^{e-1}(L^\pe)@>>>0$.
\nl
We see that

$\dim N^{e-1}(L^\pe)=\dim\ker\l-\dim\ker N^{e-1}=\dim V-\dim\ker N^{e-1}-1=c_e^N-1$
\nl
so that $c_e^{N'}=\dim N'{}^{e-1}V'$ equals $c_e^N-1$ if $c_e^N\in\ZZ'$ and equals
$c_e^N-2$ if $c_e^N\in\ZZ''$.

We write $V=W\op Y$ as in 1.9(c). Let $N_1=N|_W,N_2=N|_Y$. Let us now assume that we are 
not in the case 

($*$) $V\in\cc',\l\ne0,e=2,c_e^N\in\ZZ'$.
\nl
By 2.2(c) we have $L\sub N^{e-1}V$. Since $N^{e-1}Y=0$ we have $L\sub N^{e-1}W$ hence 
$L\sub W$. Since $\ker N^{e-1}\sub\ker\l=L^\pe$ and $Y\sub\ker N^{e-1}$ we have 
$Y\sub L^\pe$. Hence we have canonically $V'=W'\op Y$ where $W'=L^{\pe'}/L'$ (and 
$L^{\pe'}=\{x\in W;\la x,L\ra=0\}$) with both $W',Y$, $N'$-stable. Let $N'_1=N'|_{W'}$. If
$\l=0$ then $c_{e-2}^{N'_1}=c_e^{N_1}=c_e^N$ and $c_i^{N'_1}=0$ for $i\ne e-2$. If $\l\ne0$
we have either 

$c_e^{N'_1}=c_e^{N_1}-1=c_e^N-1,c_{e-2}^{N'_1}=1,c_i^{N'_1}=0$ for $i\n\{e,e-2\}$
\nl
or

$c_e^{N'_1}=c_e^{N_1}-2=c_e^N-2,c_{e-1}^{N'_1}=2,c_i^{N'_1}=0$ for $i\n\{e,e-1\}$;
\nl
the formulas for $c_i^{N'}$ in (i)-(iii) follow since 

$c_i^{N'}=c_i^{N'_1}+c_i^{N_2}=c_i^{N'_1}+c_i^N$ for $i\le e-1$.
\nl
(Those formulas hold also in the excluded case $(*)$. Indeed in that case we have 
$c_2^{N'}=c_2^N-1$ as we have seen already. For $i>2$ we have $c_i^{N'}=0$ and 
$c_1^{N'}+2c_2^{N'}=\dim V'=\dim V-2=c_1^N+2c_2^N-2$ hence $c_1^{N'}=c_1^N$.)

Assume now that $l\ne0$. We show that $\e_e^{N'}=0$. It is enough to show that if 
$x\in L^\pe$ then $\la x,N^{e-1}x\ra=0$. This follows from $L^\pe=\ker\l$. We have
$\e_{e-1}^{N'}=0$ since $e-1\in\ZZ'$. (This determines completely $\e_i^{N'}$ in the case
where $e=2$.)

Assume that $\l\ne0,c_e^N\in\ZZ',e>2$. We show that $\e_{e-2}^{N'}=1$. It is enough to show
that there exists $x\in V$ such that $N^{e-2}x\in L$ and $\la x,N^{e-3}x\ra\ne0$. (For such
$x$ we have automatically $x\in L^\pe$ since $N^{e-1}x=0$ and 
$\ker N^{e-1}\sub\ker\l=L^\pe$.) By 2.2(c) we can find $y\in V$ such that $N^{e-1}y$ is a 
basis element of $L$. We have $\la y,N^{e-1}y\ra\ne 0$. (If not we would have $y\in L^\pe$
and $L\sub N^{e-1}(L^\pe)$ contradicting 2.2(b).) Let $x=Ny$. We have 
$N^{e-2}x=N^{e-1}y\in L$ and $\la x,N^{e-3}x\ra=\la Ny,N^{e-2}y\ra=\la y,N^{e-1}y\ra\ne0$,
as required. 

We may now assume that we are not in the case $(*)$. We use the decompositions
$V=W\op Y,V'=W'\op Y$ as above. The proof of the remaining assertions on $\e_i^{N'}$ in 
(i)-(iii) is standard.

\subhead 2.5\endsubhead
Let $N\in\cm_Q$. We associate to $N$ a collection of subspaces
$(V^{\ge a})_{a\in\ZZ}=(V^{\ge a}_N)_{a\in\ZZ}$ of $V$ using induction on $D$. Let $e=e_N$.
If $N=0$ we set $V^{\ge a}=V$ for $a\le0$, $V^{\ge a}=0$ for $a\ge1$. Thus $V^{\ge a}$ are
defined when $D\le1$. We may assume that $D\ge2$, $N\ne0$ and that $V^{\ge a}$ are already
defined when $V$ is replaced by a vector space of dimension $<D$. We have $e\ge2$. Let 
$\l,L,V',N',Q'$ be as in 2.2, 2.3. Let $\r:L^\pe@>>>V'$ be the obvious map. Since $L\ne0$,
we have $\dim V'<D$. By the induction hypothesis, $V'{}^{\ge a}=V'_{N'}{}^{\ge a}$ is 
defined for any $a\in\ZZ$. If $\l=0$ we set 
$$V^{\ge a}=V\text{ for }a\le1-e;V^{\ge a}=\r\i(V'{}^{\ge a})\text{ for }a\in[2-e,e-1];
V^{\ge a}=0\text{ for }a\ge e.$$
If $\l\ne0$ we set
$$V^{\ge a}=V\text{ for }a\le-e;V^{\ge a}=\r\i(V'{}^{\ge a})\text{ for }a\in[1-e,e];
V^{\ge a}=0\text{ for }a\ge1+e.$$
This completes the definition of the subspaces $V^{\ge a}$.

From the definition it is clear that $V^{\ge a+1}\sub V^{\ge a}$ for any $a$. Thus
$(V^{\ge a})$ is a filtration of $V$. We show:

(a) {\it If $a\ge1$ we have $Q|_{V^{\ge a}}=0$ and $V^{\ge1-a}=(V^{\ge a})^\pe$.}
\nl
We use induction on $D$. For $N=0$ the result is obvious. Thus the result holds when 
$D\le1$. Now assume that $D\ge2$, $N\ne0$. Let $\l,L,V',Q',N'$ be as above. By the 
induction hypothesis we have $Q'|_{V'{}^{\ge a}_{N'}}=0$ and 
$V'{}^{\ge1-a}_{N'}=(V'{}^{\ge a}_{N'})^{\pe'}$ with ${}^{\pe'}$ relative to $V'$. It 
follows that $Q|_{V^{\ge a}}=0$ and $V^{\ge1-a}=(V^{\ge a})^\pe$ for $a\in[1,e-1]$ (if 
$\l=0$) and for $a\in[1,e]$ (if $\l\ne0$). For $a\ge e$ (with $\l=0$) and for $a\ge1+e$ 
(with $\l\ne0$) the result is again true since $Q|_{\{0\}}=0$ and $\{0\}^\pe=V$. This 
proves (a).

\subhead 2.6\endsubhead
Let $N\in\cm_Q-\{0\}$. Let $\l,L,V',N',Q',e'$ be as in 2.2, 2.3. Let 
$\r:L^\pe@>>>V'$ be the obvious map. Define $\pe'$ in terms of $V'$ in the same way as 
$\pe$ is defined in terms of $V$. If $e'\ge1$ we set $\l'=\l^{N'}_{e'}$. If $e'\ge2$ we 
define $L',V'',N'',Q'',e''$ in terms of $V',N',Q'$ in the same way as $L,V',N',Q',e'$ were
defined in terms of $V,N,Q$. If $e'\ge2$ and $e''\ge1$ we set $\l''=\l^{N''}_{e''}$.

Let $V=W\op Y$ be as in 1.9(c). For a subspace $Z$ of $Y$ let 
$Z^\vda=\{y\in Y;\la y,Z\ra=0\}$. 

In the case where $\l=0,\e'=\e-1\ge2,\l'\ne0$ we can view $L'$ as a line in $Y$ as follows.
We have $N^{e-1}V\sub W$, $Y\sub\ker N^{e-1}$ hence $V'=W'\op Y$ where 
$W'=(\ker N^{e-1}\cap W)/N^{e-1}V$. Since $W'\sub\ker N'{}^{e-2}\sub\ker\l'$ we have
$(\ker\l')^{\pe'}\sub W'{}^{\pe'}=Y$. Hence $L'\sub Y$ and $L'{}^{\pe'}=W'\op L'{}^\vda$.

In this subsection we describe explicitly $V^{\ge a}$ for certain $a$.

(i) {\it If $\l\ne0$ then $V^{\ge e}=L$.}
\nl
We have $e'\le e$. If $e'<e$ then $V'{}^{\ge e}=0$. If $e'=e$ then (using 2.4(ii),(iii)) we
have $\e_e^{N'}=0$ hence $\l'=0$ and $V'{}^{\ge e}=0$. Hence $V^{\ge e}=\r\i(0)=L$.

(ii) {\it If $\l\ne0,e'<e$ then $V^{\ge e-1}=L$.}
\nl
If $e'<e-1$ then $V'{}^{\ge e-1}=0$. If $e'=e-1$ and $\l'=0$ then $V'{}^{\ge e-1}=0$. If 
$e'=e-1$ and $\l'\ne0$ then $e',e$ are even, contradiction. Hence $V^{\ge e-1}=\r\i(0)=L$.

(iii) {\it If $\l=0$ and either $e'\le e-2$ or $e'=e-1,\l'=0$ then $V^{\ge e-1}=N^{e-1}V$.}
\nl
If $e'\le e-2$ we have $V'{}^{\ge e-1}=0$. If $e'=e-1,\l'=0$ we have again 
$V'{}^{\ge e-1}=0$. Hence $V^{\ge e-1}=\r\i(0)=L=N^{e-1}V$.

(iv) {\it If $\l\ne0,e'=e$ then $V^{\ge e-1}=N^{e-1}L^\pe+L\sub\ker N$.}
\nl
Using 2.4(ii),(iii) we have $\l'=0$. Hence $e''<e'$. If either $e''\le e'-2$ or 
$e''=e'-1,\l''=0$ then from (iii) applied to $V'$ we get $V'{}^{\ge e-1}=N'{}^{e-1}V'$ and
the result follows. If $e''=e'-1,\l''\ne0$ then $e'',e$ are even, contradiction. We have 
$NV^{\ge e-1}=0$ since $N^e=0$ and $NL=0$.

(v) {\it If $\l=0$ and either $e'\le e-3$ or $e'=e-2\ge1,\l'=0$ or $e=2,e'=0$ then 
$V^{\ge e-2}=N^{e-1}V$.}
\nl
If $e'\le e-3$ then $V'{}^{\ge e-2}=0$. If $e'=e-2\ge1,\l'=0$ then $V'{}^{\ge e-2}=0$. If 
$e=2,e'=0$ then $V'=0$ and $V'{}^{\ge e-2}=0$. Hence $V^{\ge e-2}=\r\i(0)=L=N^{e-1}V$.

(vi) {\it If $\l=0,e'=e-1,\l'\ne0$ then
$V^{\ge e-1}=\{x\in V;N^{e-1}x=0,\la x,N^{e-2}x\ra=0\}^\pe\cap q\i(0)$,}

$V^{\ge e-1}=N^{e-1}V\op L'\sub\ker N$.
\nl
By (i) for $V'$ we have 
$V'{}^{\ge e-1}=L'=\{x\in V';\la x,N^{e-2}x\ra=0\}^{\pe'}\cap q'{}\i(0)$. Hence 

$V^{\ge e-1}=\ker N^{e-1}\cap\{x\in\ker N^{e-1};\la x,N^{e-2}x\ra=0\}^\pe\cap q\i(0)$.
\nl
In the last equality of (vi) we regard $N^{e-1}V$ as a subspace of $W$ and $L'$ as a 
subspace of $Y$ (as earlier in this subsection).

(vii) {\it If $\l=0,e'=e-2\ge1,\l'\ne0$ we have 
$V^{\ge e-2}=\{x\in V;N^{e-1}x=0,\la x,N^{e-3}x\ra=0\}^\pe\sub\ker N$.}
\nl
By (i) for $V'$ we have 
$V'{}^{\ge e-2}=\{x'\in V';\la x',N'{}^{e-3}x'\ra0\}^{\pe'}\cap Q'{}\i(0)$. Hence 

$V^{\ge e-2}=\ker N^{e-1}\cap\{x\in\ker N^{e-1};\la x,N^{e-3}x\ra=0\}^\pe\cap Q\i(0)$.
\nl
We have $NV\sub\{x\in\ker N^{e-1};\la x,N^{e-3}x\ra=0\}$ since $\l=0$. Taking $\pe$ we 
obtain $\{x\in\ker N^{e-1};\la x,N^{e-3}x\ra=0\}^\pe\sub\ker N$. Since 
$\ker N\sub\ker N^{e-1}$ we see that $V^{\ge e-2}$ is as required.

(viii) {\it If $\l=0,e'=e-1,\l'\ne0,e''\le e'-1$, we have 
$V^{\ge e-2}=\{x\in V;N^{e-1}x=0,\la x,N^{e-2}x\ra=0\}^\pe\cap Q\i(0)\sub\ker N$.}
\nl
By (ii) for $V'$ we have 
$V'{}^{\ge e-2}=\{x'\in V';\la x',N'{}^{e-2}x'\ra0\}^\pe\cap Q'{}\i(0)$. Hence

$V^{\ge e-2}=\ker N^{e-1}\cap\{x\in\ker N^{e-1};\la x,N^{e-2}x\ra=0\}^\pe\cap Q\i(0)$.
\nl
We have $NV\sub\{x\in\ker N^{e-1};\la x,N^{e-2}x\ra=0\}$. Taking $\pe$ we obtain
$\{x\in\ker N^{e-1};\la x,N^{e-2}x\ra=0\}^\pe\sub\ker N$. Since $\ker N\sub\ker N^{e-1}$ we
see that $V^{\ge e-2}$ is as required.

(ix) {\it If $\l=0,e'=e-1,\l'=0$ and either $e'=1$ or $e\ge3,e''\le e'-2$ or 
$e\ge3,e''=e'-1,\l''=0$ then $V^{\ge e-2}=N^{e-2}(\ker N^{e-1})+N^{e-1}V\sub\ker N$.}
\nl
By (iii) for $V'$ we have $V'{}^{\ge e-2}=N'{}^{e-2}V'$. Hence $V^{\ge e-2}$ is as 
required.

(x) {\it If $\l=0,e'=e-1,\l'\ne0,e''=e'$, then 
$V^{\ge e-2}=N^{e-1}V\op(N^{e-2}L'{}^\vda+L')\sub\ker N$.}
\nl
Here we regard $N^{e-1}V$ as a subspace of $W$ and $L'$ as a subspace of $Y$ (as earlier in
this subsection). Then $N^{e-2}L'{}^\vda+L'$ is a subspace of $Y$. By (iv) for $V'$ we have
$V'{}^{\ge e-2}=N'{}^{e-2}L'{}^{\pe'}+L'$. Since $N'{}^{e-2}W'=0$ we have 
$V'{}^{\ge e-2}=N'{}^{e-2}L'{}^\vda+L'\sub Y$ and 
$V^{\ge e-2}=N^{e-1}V\op(N^{e-2}L'{}^\vda+L')$. By (iv) for $V'$ we have 
$N(N^{e-2}L'{}^\vda+L')=0$. Hence $NV^{\ge e-2}=0$.

(xi) {\it If $\l=0,e'=e-1,e\ge3,\l'=0,e''=e'-1,\l''\ne0$ then
$V^{\ge e-2}=N^{e-1}V\op(\{b\in Y;N^{e-2}b=0,\la b,N^{e-3}b\ra=0\}^\vda\cap q\i(0))
\sub\ker N$.}
\nl
By (vi) for $V'$ we have $V'{}^{\ge e-2}=U'{}^{\pe'}\cap Q'{}\i(0)$ where 
$U'=\{x\in V';N'{}^{e-2}x=0,\la x,N'{}^{e-3}x'=0\}$. We write $V'=W'\op Y$ as in the proof
of (x). We write the condition that $x=a+b$ with $a\in W',b\in Y$ is in $U'$ in terms of
$a,b$. Note that $N'{}^{e-2}a=0$ and $\la a,N'{}^{e-3}a\ra=0$. (The last equality follows 
from $\ta\in W,N^{e-1}\ta=0\imp\la\ta,N^{e-3}\ta\ra=0$. Indeed we have $\ta=Nc$ with
$c\in W$ and $\la\ta,N^{e-3}\ta\ra=\la Nc,N^{e-3}Nc\ra=\la c,N^{e-1}c\ra=0$ since $\l=0$.)
We see that $U'=W'\op\{b\in Y;N^{e-2}b=0,\la b,N^{e-3}b\ra=0\}$ and

$V'{}^{\ge e-2}=\{b\in Y;N^{e-2}b=0,\la b,N^{e-3}b\ra=0\}^\vda\cap Q'{}\i(0)$,

$V^{\ge e-2}=N^{e-1}V\op(\{b\in Y;N^{e-2}b=0,\la b,N^{e-3}b\ra=0\}^\vda\cap Q\i(0))$.
\nl
By (vi) for $V'$ we have $N(\{b\in Y;N^{e-2}b=0,\la b,N^{e-3}b\ra=0\}^\vda\cap Q\i(0))=0$.
Hence $NV^{\ge e-2}=0$.

\subhead 2.7\endsubhead
Let $N\in\cm_Q,e=e_N$. Let $V^{\ge a}$ be as in 2.5. We show:

(a) {\it $NV^{\ge a}\sub V^{\ge a+2}$ for any $a\in\ZZ$.}
\nl
When $N=0$ the result is obvious. Now assume that $N\ne0$. Then $e\ge2$. Let 
$\l,L,V',N',Q'$ be as in 2.2, 2.3. We may assume that (a) holds when $V,N$ are replaced by
$V',N'$. We may assume that (a) holds when $V,N$ are replaced by $V',N'$.

Assume first that $\l=0$. If $a\ge e$ then $V^{\ge a}=0$ and (a) is obvious. If 
$a\in\{e-2,e-1\}$ then $NV^{\ge a}=0$ by 2.6 and (a) holds. Assume now that $a=-e$ or that 
$a=1-e,e\ge3$. To prove (a) in this case it is enough to show that $NV\sub V^{\ge a+2}$ 
that is (using 2.5(a)) $NV\sub(V^{\ge-1-a})^\pe$ or that $N^\da(V^{\ge-1-a})\in R$. This 
follows from $N(V^{\ge-1-a})=0$ which has been noted earlier. If $a=-1,e=2$ we have 
$NV=V^{\ge a+2}$ and (a) holds. If $a\le-1-e$ then $V^{\ge a+2}=V$ and (a) is obvious. If 
$2-e\le a\le e-3$ then $V^{\ge a}=\r\i(V'{}^{\ge a})$, $V^{\ge a+2}=\r\i(V'{}^{\ge a+2})$ 
(notation of 2.5). Since $N'V'{}^{\ge a}\sub V'{}^{\ge a+2}$ we see that (a) holds.

Assume next that $\l\ne0$. If $a\ge e+1$ then $V^{\ge a}=0$ and (a) is obvious. If 
$a\in\{e-1,e\}$ then $NV^{\ge a}=0$ by 2.6 and (a) holds. Assume now that 
$a\in\{-e,-1-e\}$. To prove (a) in this case it is enough to show that $NV\sub V^{\ge a+2}$
that is (using 2.5(a)) $NV\sub(V^{\ge-1-a})^\pe$ or that $N^\da(V^{\ge-1-a})\in R$. This 
follows from $N(V^{\ge-1-a})=0$ which has been noted earlier. If $a\le-2-e$ then 
$V^{\ge a+2}=V$ and (a) is obvious. If $1-e\le a\le e-2$ then 
$V^{\ge a}=\r\i(V'{}^{\ge a})$, $V^{\ge a+2}=\r\i(V'{}^{\ge a+2})$ (notation of 2.5). Since
$N'V'{}^{\ge a}\sub V'{}^{\ge a+2}$ we see that (a) holds. This proves (a).

For any $a\in\ZZ$ we set $\bV^a=V^{\ge a}/V^{\ge a+1}$. From (a) we see that $N$ induces a
linear map $\bN:\bV^a@>>>\bV^{a+2}$.

\subhead 2.8\endsubhead
In the setup of 2.6 taking $\pe$ in 2.6(i),(iii),(vi) we obtain

{\it If $\l\ne0$ then $V^{\ge 1-e}=L^\pe$.}

{\it If $\l=0$ and either $e'\le e-2$ or $e'=e-1,\l'=0$ then 
$V^{\ge 2-e}=\ker N^{e-1}$.}

{\it If $\l=0,e'=e-1,\l'\ne0$ then $V^{\ge 2-e}=(\ker N^{e-1}\cap W)\op L'{}^\vda$.}

\subhead 2.9\endsubhead
For any $a\ge0$ we set $K_a=\ker\bN^a:\bV^{-a}@>>>\bV^a$. For any $a\in\ZZ''_{\ge2}$ we 
define a quadratic form $Q_a:\bV^{-a}@>>>\kk$ by $Q_a(x)=Q_0(\bN^{a/2}x)$. (Note that 
$\bN^{a/2}x\in\bV^0$.) If $\dx$ is a representative of $x$ in $V^{\ge-a}$ we have 

$Q_a(x)=Q(N^{a/2}\dx)=\la N^{a/2-1}\dx,N^{a/2}\dx\ra=
\la\dx,N^{a-1}\dx+cN^a\dx+\do\ra$
\nl
where $c\in\kk$. If $f>0$ we have $\la x,N^{a+f}x\ra=0$ since 
$\la V^{\ge-a},V^{\ge a+2f}\ra=0$. Moreover,
$$\align&\la x,N^ax\ra=\la N^{a/2}x+c'N^{a/2+1}x+\do,N^{a/2}x\ra\\&=
\la c'N^{a/2+1}x+\do,N^{a/2}x\ra=c'\la x,N^{a+1}x\ra=0\endalign$$
where $c'\in\kk$. Hence
$$Q_a(x)=\la\dx,N^{a-1}\dx\ra.$$
Let $\la,\ra_a$ be the symplectic form on $\bV^{-a}$ associated to $Q_a$. Thus 
$$\la x,x'\ra_a=\la\bN^{a/2}x,\bN^{a/2}x'\ra_0$$
for $x,x'\in\bV^{-a}$. If $\dx,\dx'$ are representatives of $x,x'$ in $V^{\ge-a}$ we have
$$\la x,x'\ra_a=\la N^{a/2}\dx,N^{a/2}\dx'\ra=\la\dx,N^a\dx'+cN^{a+1}\dx'+\do\ra
=\la\dx,N^a\dx'\ra=\la x,\bN^a x'\ra.$$
(We use that $\la V^{\ge-a},V^{\ge a+2}\ra=0$.) Let $R_a$ be the radical of $\la,\ra_a$. If
$x'\in R_a$ then $\la x,\bN^ax'\ra=0$ for all $x\in\bV^{-a}$ hence $\bN^ax'=0$. Thus 
$R_a=K_a$. We show:

(a) {\it If $a\in\ZZ'$ we have $K_a=0$. If $a\in\ZZ''_{\ge0}0$ then $Q_a$ is nondegenerate;
hence $\dim K_a\in\{0,1\}$.}
\nl
If $N=0$ we have $V^{\ge-a}=V_{\ge1-a}=V$ hence $\bV^{-a}=0$ so that $K_a=0$ as required.
Now assume that $N\ne0$ so that $e\ge2$. Let $\l,L,V',N',Q'$ be as in 2.2, 2.3. Let 
$\bV'{}^a,\bN'{}^a,Q'_a$ be the analogues of $\bV^a,\bN^a,Q_a$ for $V',N'$ instead of 
$V,N$. We may assume that the analogue of (a) holds when $V,N$ is replaced by $V',N'$. If 
$\l\ne0,a>e$ or if $\l=0,a\ge e$ we have $\bV^{-a}=0$ hence $K_a=0$ as required. If 
$\l\ne0,a=e$ we have $\bV^{-e}=V/L^\pe$, $\bV^e=L$, $bN^e=0$ hence $K_a=V/L^\pe$. We must 
show that $Q_a$ is not identically zero. It is enough to show that 
$\dx\m\la\dx,N^{a-1}\dx\ra$ is not identically zero on $V$; this holds since $\l\ne0$. If 
$\l=0,a=e-1,\e_{e-1}^{N'}=0$ then $\bV^{1-e}=V/\ker N^{e-1}$, $\bV^{e-1}=N^{e-1}V$ and 
$\bN^a:\bV^{-a}@>>>\bV^a$ is an isomorphism; hence $K_a=0$, as required. If 
$\l=0,a=e-1,\e_{e-1}^{N'}=1$ then with notation of 2.6(vi), 2.8 we have
$\bV^{1-e}=W/(\ker N^{e-1}|W)\op Y/L'{}^\vda$, $\bV^{e-1}=N^{e-1}V\op L'$; now $\bN^{e-1}$ 
restricts to an isomorphism $W/(\ker N^{e-1}|W)@>>>N^{e-1}V$ and to the zero map 
$Y/L'{}^\vda@>>>L'$ (since $N^{e-1}Y=0$). Hence $K_a=Y/L'{}^\vda$. We must show that $Q_a$
is not identically $0$ on $Y/L'{}^\vda$ or that $x\m\la x,N^{e-2}x\ra$ is not identically 
$0$ on $Y$. But this follows from $\e_{e-1}^{N'}=1$. If $\l\ne0,a\in[1,e-1]$ or if $\l=0$,
$a\in[1,e-2]$, we have $\bV^{-a}=\bV'{}^{-a}$, $\bV^a=\bV'{}^a$. We can identify
$\bN'{}^a=\bN^a$ and (if $a\in\ZZ''$) $Q_a$ with $Q'_a$. Hence the result follows from the
induction hypothesis. If $a=0$ the result is obvious. This completes the proof.

We show:

(b) {\it Assume that $a\in\ZZ''_{\ge2}$. We set $\ph_a=|\{b\in\ZZ'';c_b^N\in\ZZ',b>a\}|$,
$\x_a=\dim K_a$. We have $c_a^N\in\ZZ'\imp\x_a=1$; $c_a^N\in\ZZ'',\ph_a\in\ZZ'\imp\x_a=1$;
$c_a^N\in\ZZ'',\ph_a\in\ZZ''\imp\x_a=\e_a^N$.}
\nl
If $N=0$ the result is obvious. Now assume that $N\ne0$ so that $e\ge2$. Let $\l,V',N',Q'$
be as in 2.2, 2.3. Let $\x'_a,\ph'_a$ be the analogues of $\x_a,\ph_a$ for $V',N'$ instead
of $V,N$. We may assume that the analogue of (b) holds when $V,N$ is replaced by $V',N'$. 
If $\l\ne0,a>e$ or if $\l=0,a\ge e$ we have (by the proof of (a)) $\x_a=0$, as required. If
$\l\ne0,a=e$ we have (by the proof of (a)) $\x_a=1$ as required. If $\l=0,a=e-1$ then (by 
the proof of (a)) $\x_a=\e_{e-1}^{N'}=\e_{e-1}^N$ as required. In the remainder of the 
proof we assume that either $\l\ne0,a\in[1,e-1]$ or $\l=0$, $a\in[1,e-2]$. Then (by the 
proof of (a)), $\x_a=\x'_a$. Using the induction hypothesis we see that if 
$c_a^{N'}\in\ZZ'',\ph'_a\in\ZZ''$ then $\x_a=\e'_a{}^N$; otherwise, $\x_a=1$.

Assume first that $\l\ne0$. Then $e\in\ZZ''$ hence $a\le e-2$.

If $a=e-2$, $c_{e-2}^N\in\ZZ'$ and $c_e^N\in\ZZ'$ then $c_{e-2}^{N'}\in\ZZ''$, 
$c_e^{N'}\in\ZZ''$ and $\e^{N'}_{e-2}=1$ so that $\x_{e-2}=1$. 

If $a=e-2$, $c_{e-2}^N\in\ZZ'$ and $c_e^N\in\ZZ''$ then $c_{e-2}^{N'}\in\ZZ'$, 
$c_e^{N'}\in\ZZ''$ so that $\x_{e-2}=1$. 

If $a=e-2$, $c_{e-2}^N\in\ZZ''$ and $c_e^N\in\ZZ'$ then $c_{e-2}^{N'}\in\ZZ'$, 
$c_e^{N'}\in\ZZ''$ and $\e^{N'}_{e-2}=1$ so that $\x_{e-2}=1$. Also, $\ph_{e-2}=1$.

If $a=e-2$, $c_{e-2}^N\in\ZZ''$ and $c_e^N\in\ZZ''$ then $c_{e-2}^{N'}\in\ZZ''$, 
$c_e^{N'}\in\ZZ''$, $x_{e-2}=\e_{e-2}^{N'}=\e_{e-2}^N$. 

If $a\le e-4$ we have $c_a^{N'}=c_a^N$. If these are odd then $\x_a=1$. If these are even 
then $\ph_a=\ph'_a\mod 2$ (if $c_e^N\in\ZZ'$, $c_{e-2}^N\in\ZZ'$ then $\ph_a=\ph'_a+2$; 
otherwise, $\ph_a=\ph'_a$).  

Assume next that $\l=0$. 

If $a=e-2$ then $c_e^N\in\ZZ''$ and $c_a^{N'}=c_a^N\mod2$. If $c_a^{N'},c_a^N$ are odd we 
have $\x_a=1$. If $c_a^{N'},c_a^N$ are even we have $\x_a=\e^{N'}_a=\e^N_a$. 

If $a\le e-3$ we have $c_a^{N'}=c_a^N$. If these are both odd then $\x_a=1$; if these are 
both even then $\ph_a=\ph'_a$. 

We see that (b) holds.

We show:

(c) {\it If $a\in ZZ'_{\ge1}$ then the bilinear pairing $\bV^{-a}\T\bV^{-a}@>>>\kk$, 
$x,y\m\la x,\bN^a y\ra$ is symplectic; it is nondegenerate by (a). Hence 
$\bV^{-a}\in\cc''$.}
\nl
If $N=0$ we have $\bV^{-a}=0$ if $a\ne0$ and (c) is obvious. Assume that $N\ne0$. Then 
$e\ge2$. Let $\l,V',N'$ be as in 2.2, 2.3. We may assume that (c) holds when $V,N$ are
replaced by $V',N'$. If $\l\ne0,a>e$ or if $\l=0,a\ge e$ we have $\bV^{-a}=0$ and (c) is 
obvious. The case $\l\ne0,a=e$ cannot occur since $e\in\ZZ''$. If $\l=0,a=e-1$ then 
$e-1\in\ZZ'$ and $\e_{e-1}^{N'}=0$; we have $\bV^{1-e}=V/\ker N^{e-1}$ and 
$\la x,N^{e-1}x\ra=0$ for all $x$. Hence (c) holds. If $\l\ne0,a\in[1,e-1]$ or if 
$\l=0,a\in[1,e-2]$, we have $\bV^{-a}=\bV'{}^{-a}$. Hence the result follows the induction
hypothesis. The case $a=0$ does not arise. This completes the proof.

\subhead 2.10\endsubhead
Let $X^*=(X^{\ge a})_{a\in\ZZ}$ be a $Q$-filtration of $V$ (see 1.4). Let 
$N\in\cm_Q,e=e_N$. We say that $X^*$ is {\it $N$-adapted} if conditions (i)-(iii) below 
hold.

(i) $NX^{\ge a}\sub X^{\ge a+2}$ for any $a$;
\nl
For any $a\ge0$ let $\ck_a$ be the kernel of the map $\nu^a:\gr^{-a}X^*@>>>\gr^aX^*$ 
induced by $N^a$. We state conditions (ii),(iii).

(ii) for any $a\in ZZ'_{\ge1}$ we have $\ck_a=0$;

(iii) for any $a\in\ZZ''_{\ge2}$ we have $\ck_a=0$ or $\dim\ck_a=1$; 
in the latter case the map $\ck_a@>>>\kk$, $x\m\la\dx,N^{a-1}\dx\ra$ is a bijection.
\nl
(For $x\in\gr^aX^*$ we denote by $\dx$ a representative of $x$ in $X^{\ge a}$.) Note that
$(V^{\ge a})_{a\in\ZZ}$ is $N$-adapted where $V^{\ge a}=V^{\ge a}_N$. We show:

(a) {\it If $(X^{\ge a})_{a\in\ZZ}$ is an $N$-adapted filtration of $V$ then
$X^{\ge a}=V^{\ge a}$ for any $a$.}
\nl
If $V=0$ the result is obvious. Now assume that $V\ne0$. If $N=0$ and $a\ge1$ then
$\nu^a:\gr^{-a}X^*@>>>\gr^aX^*$ is $0$ and $\ck_a=0$. Hence $X^{\ge a}=X^{\ge a+1}$ and 
$X^{\ge-a }=X^{\ge-a+1}$ so that we have $X^{\ge1}=X^{\ge2}=\do=0$ and 
$X_{\ge0}=X_{\ge-1}=\do=V$ as required. We now assume that $N\ne0$ hence $e\ge2$. Let
$\l,L,V',N',Q',\r,\e'$ be as in 2.2, 2.3. If $\e'\ge1$ let $\l'$ be as in 2.6. If $\e'\ge2$
let $Y,L',L'{}^\vda$ be as in 2.6. We may assume that (a) holds when $V,N$ are replaced by
$V',N'$.

If $a\ge e+1$ or if $a=e,\l=0$, then $\nu^a:\gr^{-a}X^*@>>>\gr^aX^*$ is $0$ and $\ck_a=0$. 
Hence $X^{\ge a}=X^{\ge a+1}$ and $X^{\ge-a }=X^{\ge-a+1}$ so that if $\l=0$ we have

$X^{\ge e}=X^{\ge e+1}=\do=0$ and $X_{\ge1-e}=X_{\ge-e}=\do=V$;
\nl
if $\l\ne0$ we have

$X^{\ge e+1}=X^{\ge e+2}=\do=0$ and $X_{\ge-e}=X_{\ge-1-e}=\do=V$.
\nl
Hence 
$X^{\ge a}=V^{\ge a}$ if $a\ge e+1$ or if $a=e,\l=0$ or if $a\le-e$ or if $a=1-e,\l=0$. 

Assume that $\l\ne0$. Then $\ck_e$ is the kernel of 
$0=\nu^e:V/X^{\ge 1-e}@>>>X^{\ge e}$ that is $\ck_e=V/X^{\ge 1-e}$. Also $\ck_e@>>>\kk$, 
$x\m\la\dx,N^{e-1}\dx\ra$ is not identically zero hence $\dim V/X^{\ge 1-e}=1$. Since for 
$x\in\gr^{1-e}X^*$ we have $\dx,N^{e-1}\dx\ra=0$ we see that $X^{\ge 1-e}\sub\ker\l=L^\pe$.
Since $\dim X^{\ge 1-e}=\dim L^\pe$ we see that $X^{\ge 1-e}=L^\pe$. We have 
$X^{\ge e}=(X^{\ge 1-e})^\pe\cap Q\i(0)=L$ as required.

Next we assume that $\l=0$. Then $\ck_{e-1}$, the kernel of 
$\nu^a:V/X^{\ge 2-e}@>>>X^{\ge e-1}$, is
$$\align&\{x\in V/X^{2-e};N^{e-1}\dx\in X^{\ge e}=0\}=\{x\in V/X^{2-e}; N^{e-1}\dx=0\}
\\&=\ker N^{e-1}/X^{\ge2-e}.\endalign$$
If $\e_{e-1}^N=0$ then $\dx,N^{e-2}\dx\ra=0$ for $x$ in this kernel that is for
$x\in\ck_{e-1}$. Hence in this case we have $\ck_{e-1}=0$. Thus $N^{e-1}$ induces an 
isomorphism $V/X^{\ge 2-e}@>>>X^{\ge e-1}$ so that 
$X^{\ge e-1}=N^{e-1}V,X^{\ge2-e}=\ker N^{e-1}$, as required. We have 
$X^{\ge 2-e}=(X^{\ge e-1})^\pe=(V^{\ge e-1})^\pe=V^{\ge 2-e}$ as required.

Now assume that $\l=0,\e_{e-1}^N=1$. In this case we have an isomorphism 
$\ker N^{e-1}/X^{2-e}@>>>\kk$ induced by $\l_{e-1}^N$ that is we have 

$X^{\ge2-e}=\{x\in\ker N^{e-1};\l_{e-1}^Nx=0\}$.
\nl
This is the same as $(\ker N^{e-1}\cap W)\op L'{}^\vda=V^{\ge2-e}$, see 2.8. Taking $\pe$ 
in $X^{\ge2-e}=V^{\ge2-e}$ and intersecting with $Q\i(0)$ we obtain 
$X^{\ge e-1}=V^{\ge e-1}$.

If $\l=0,a\in[2-e,e-1]$ we have $N^{e-1}V\sub X^{\ge a}\sub\ker N^{e-1}$ and we denote by 
$X'{}^{\ge a}$ the image of $X^{\ge a}$ under 
$\r:\ker N^{e-1}@>>>V'=\ker N^{e-1}/N^{e-1}V$. For $a\le1-e$ we set $X'{}^{\ge a}=V'$ and 
for $a\ge e$ we set $X'{}^{\ge a}=0$. Now $(X'{}^{\ge a})_{a\in\ZZ}$ is an $N'$-adapted
filtration of $V'$. (We must only show the analogue of (ii),(iii) for $N'$ with $a=e-1$. We
have $X'{}^{\ge1-e}/X'{}^{\ge2-e}=\ker N^{e-1}/X^{\ge 2-e}$; this is $Y/L'{}^\vda$ if 
$\e'=\e-1,\l'\ne0$ and $0$ otherwise. We have
$X'{}^{\ge e-1}/X'{}^{\ge e}=X^{\ge e-1}/N^{e-1}V$. This is $L'$ if $\e'=\e-1,\l'\ne0$ and
$0$ otherwise. Hence (ii),(iii) are obvious in this case.) By the induction hypothesis we 
have $X'{}^{\ge a}=V'{}^{\ge a}_N$ for all $a$. Taking inverse image under $\r$ we see that
for $a\in[2-e,e-1]$ we have $X^{\ge a}=V^{\ge a}$.

If $\l\ne0,a\in[1-e,e]$ we have $L\sub X^{\ge a}\sub L^\pe$ and we denote by $X'{}^{\ge a}$
the image of $X^{\ge a}$ under $\r:L^\pe@>>>V'=L^\pe/L$. For $a\le-e$ we set
$X'{}^{\ge a}=V'$ and for $a\ge e+1$ we set $X'{}^{\ge a}=0$. Now 
$(X'{}^{\ge a})_{a\in\ZZ}$ is an $N'$-adapted filtration of $V'$. (We have 
$X'{}^{\ge-e}/X'{}^{\ge1-e}=0$, $X'{}^{\ge e}/X'{}^{\ge e+1}=0$.) By the induction 
hypothesis we have $X'{}^{\ge a}=V'{}^{\ge a}_N$ for all $a$. Taking inverse image under 
$\r$ we see that for $a\in[1-e,e]$ we have $X^{\ge a}=V^{\ge a}$. This completes the proof
of (a).

\subhead 2.11\endsubhead
Assume that $V\in\cc''$. If $g\in O_Q$ and $S\in\cj_Q$, we set 
$\d_g=\dim(S/(S\cap g(S)))\mod2$. It is known that $\d_g$ is independent of the choice of 
$S$ and that $\d_g=0$ if and only if $g\in SO_Q$.

Let $N\in\tcm_Q$. We show that 

(a) $\d_{1+N}=\dim\ker N\mod2$.
\nl
If $N=0$ this is clear. Assume now that $N\ne0$. Let $e=e_N,L=L_N$. We have $e\ge2$. Let
$\l,L$ be as in 2.2. We may assume that $V\in\cc''$. Assume first that $e\ge3$. As in the 
proof of 2.2(e) we have $Q|_L=0$. We set $V'=L^\pe/L$. The nondegenerate quadratic form
$Q':V'@>>>\kk$ can be defined as in 2.3. The nilpotent endomorphism $N':V'@>>>V'$ induced
by $N$ belongs to $\tcm_{Q'}$. As in the proof of 2.3(a) we see that 
$\dim\ker N'=\dim\ker N\mod2$. We have $\d_{1+N}=\d_{1+N'}$. Since the result may be 
assumed to hold for $N'$ we see that (a) holds. We now assume that $e=2$ that is 
$N^2=0,N\ne0$. 

Assume that $\l\ne0$. We can find $x\in V$ such that $\la x,Nx\ra\ne0$. Then $x,Nx$ span a
two-dimensional $N$-stable subspace $P$ of $V$ on which $\la,\ra$ is non-degenerate. Let 
$V'=P^\pe$ and let $Q'=Q|_{V'}$. Then $V=P\op V'$, $Q'$ is non-degenerate and $N$ restricts
to a nilpotent map $N':V'@>>>V'$. Note that $\dim\ker N=\dim\ker N'+1$ and 
$\d_{1+N}=\d_{1+N'}+1\mod2$. Since the result may be assumed to hold for $N'$ we see that 
(a) holds.

Assume that $\l=0$. We write $V=W\op Y$ as in 1.9(c). Then $W,Y$ are $N$-stable
nondegenerate even dimensional subspaces of $V$ with $\la W,Y\ra=0$; moreover $NY=0$. Hence
$\dim\ker N=\dim\ker(N|_W)+\dim Y$, $\d_{1+N}=\d_{1+N|_W}$. If $Y\ne0$ we may assume that the 
result holds for $N|_W$; we see that (a) holds. Thus we may assume that $Y=0$. We have
$V=E\op NE$ with $E$ as in 1.9(a). Note that $\dim NE=D/2$, see 1.9(b). Since $\l=0$ we 
have $Q|_{NE}=0$. Clearly $NE$ is $1+N$ stable. We see that $\d_{1+N}=0$. The nondegenerate
symmetric bilinear form on $E$ described in 1.9(a) is symplectic since $\l=0$. Hence 
$\dim E\in\ZZ''$ and $\dim NE\in\ZZ''$. We see that $\dim\ker N\in\ZZ''$. Thus (a) holds.

\mpb

From (a) we deduce that for $N\in\tcm_Q$ we have $N\in\cm_Q$ if and only if $1+N\in SO_Q$.

\subhead 2.12\endsubhead
We prove Theorem 1.7 (with $p=2$) in the form 1.7(a). Let $T\in SO_Q$ be unipotent. Then
$N=T-1\in O_Q$ and by 2.11, we have $N=T-1\in\cm_Q$. In 2.5 we have attached to $N$ a 
$Q$-filtration $X^*=(V^{\ge a}_N)$ of $V$. In 2.7, 2.9 we have shown that 
$N\in E^{\ge2}_*X^*$. In 2.10 we have shown that the last property determines $X^*$ 
uniquely. Thus 1.7(a) is established.

\head 3. The case $p\ne2$\endhead
\subhead 3.1\endsubhead
In this section we assume that $p\ne2$. In this case we have $R=0$.

Let $N\in\tcm_Q-\{0\},e=e_N$. We have $e\ge2$. Let $L=L_N=N^{e-1}V$, a subspace of $V$. By 
1.1(b) we have $L^\pe=\ker N^{e-1}$. Since $2e-2\ge e$ we have $N^{e-1}V\sub\ker N^{e-1}$ 
hence:

(a) $L\sub L^\pe$.
\nl
Clearly,

(b) $NL=0,NV\sub L^\pe$.
\nl
We show:

(c) $Q|_L=0$.
\nl
Let $v\in V$. We have $q(N^{e-1}v)=\la N^{e-2}v,N^{e-1}v\ra=\pm\la v,N^{2e-3}v\ra$. If
$e\ge3$ this is $0$ since $2e-3\ge e$. If $e=2$ this is $0$ since 
$\la v,Nv\ra=-\la Nv,v\ra$ so that $2\la v,Nv\ra=0$ and $\la v,Nv\ra=0$.

\subhead 3.2\endsubhead
Let $N\in\tcm_Q-\{0\},e=e_N,L=L_N$. We have $e\ge2$. We set $V'=L^\pe/L$, see 3.1(a). From
3.1(b) we see that $N$ induces a (nilpotent) endomorphism of $V'$, denoted by $N'$. Let 
$e'=e_{N'}$. We have $e'\le e-1$. Define a quadratic form $Q':V'@>>>\kk$ by 
$Q'(x')=Q(x)$ where $x$ is a representative of $x'\in V'$ in $L^\pe$. (To see that $Q'$ is
well defined we use 3.1(c).) The symmetric bilinear form associated to $Q'$ is 
$\la x',y'\ra'=\la x,y\ra$ where $x',y'\in V'$ and $x,y$ are representatives of $x',y'$ in 
$L^\pe$. Its radical is $\{x\in L^\pe;\la x,L^\pe\ra=0\}/L=0$. It follows that $Q'$ is 
nondegenerate. We have $N'\in\tcm_{Q'}$. 

\subhead 3.3\endsubhead
Let $N\in\tcm_Q,e=e_N$. We associate to $N$ a collection of subspaces

$(V^{\ge a})_{a\in\ZZ}=(V^{\ge a}_N)_{a\in\ZZ}$
\nl
of $V$ using induction on $D$. If $N=0$ we set $V^{\ge a}=V$ for $a\le0$, $V^{\ge a}=0$ for
$a\ge1$. Thus $V^{\ge a}$ are defined when $D\le1$. We may assume that $D\ge2$, $N\ne0$ and
that $V^{\ge a}$ are already defined when $V$ is replaced by a vector space of dimension 
$<D$. We have $e\ge2$. Let $L,V',N',Q'$ be as in 3.1, 3.2. Let $\r:L^\pe@>>>V'$ be the 
obvious map. Since $L\ne0$, we have $\dim V'<D$. By the induction hypothesis, 
$V'{}^{\ge a}=V'_{N'}{}^{\ge a}$ is defined for any $a\in\ZZ$. We set 
$$V^{\ge a}=V\text{ for }a\le1-e;V^{\ge a}=\r\i(V'{}^{\ge a})\text{ for }a\in[2-e,e-1];
V^{\ge a}=0\text{ for }a\ge e.$$
This completes the definition of the subspaces $V^{\ge a}$.

From the definition it is clear that $(V^{\ge a})$ is a filtration of $V$. We show:

(a) {\it If $a\ge1$ we have $Q|_{V^{\ge a}}=0$ and $V^{\ge1-a}=(V^{\ge a})^\pe$.}
\nl
We use induction on $D$. For $N=0$ the result is obvious. Thus the result holds when 
$D\le1$. Now assume that $D\ge2$, $N\ne0$. Let $L,V',Q',N'$ be as above. By the induction 
hypothesis we have $Q'|_{V'{}^{\ge a}_{N'}}=0$ and 
$V'{}^{\ge1-a}_{N'}=(V'{}^{\ge a}_{N'})^{\pe'}$ with ${}^{\pe'}$ relative to $V'$. It 
follows that $Q|_{V^{\ge a}}=0$ and $V^{\ge1-a}=(V^{\ge a})^\pe$ for $a\in[1,e-1]$. For 
$a\ge e$ the result is again true since $Q|_{\{0\}}=0$ and $\{0\}^\pe=V$. This proves (a). 

We see that $V^*=(V^{\ge a})$ is a $Q$-filtration of $V$. Clearly, $V^{\ge a}$ is the same
as the subspace $V^N_{\ge a}$ defined in terms of the nilpotent endomorphism $N:V@>>>V$ 
(without reference to $Q$) in \cite{\LU, 2.3, 2.4}. It follows that 
$NV^{\ge a}\sub V^{\ge a+2}$ for any $a\in\ZZ$. For any $a\in\ZZ$ we set 
$\bV^a=V^{\ge a}/V^{\ge a+1}$. We see that $N$ induces a linear map 
$\bN:\bV^a@>>>\bV^{a+2}$. From \cite{\LU, 2.3} we see that for any $a\ge0$, 
$\bN^a:\bV^{-a}@>>>\bV^{-a}$ is an isomorphism. It follows that $N\in E^{\ge2}_*V^*$. 
Conversely if $X^*$ is a $Q$-filtration of $V$ such that $N\in E^{\ge2}_*X^*$ we see that 
for any $a\ge0$ the kernel $\ck_a$ of the map $\gr^{-a}X^*@>>>\gr^aX^*$ induced by $N^a$ is
the radical of the symmetric bilinear form attached to a nondegenerate quadratic form on 
$\gr^{-a}X^*$; since $p\ne2$ it follows that $\ck_a=0$. Hence the map 
$\gr^{-a}X^*@>>>\gr^aX^*$ induced by $N^a$ is an isomorphism. Using \cite{\LU, 2.4} it 
follows that $X^*=V^*$. Thus 1.7(a) holds.

\head 4. On unipotent conjugacy classes in $SO_Q$ ($p=2$)\endhead
\subhead 4.1\endsubhead
In this section we assume that $\kk$ is algebraically closed. Moreover in this and the next
subsection we assume that $p=2$. Assume that $D\ge2$. Let $\ph\in\tcf_1^D$. (See 1.6.) Thus
$\ph=(f_a)$ where $f_0>0$. Let 

$\fX_\ph=\{i\in\ZZ'_{\ge1};i=f_a\text{ for some }a\in\ZZ''_{\ge0}\}$. 
\nl
this is a finite set. Let $\ce_\ph$ be the set of all subsets of $\fX_\ph\T\fX_\ph$ which
are equivalence relations on $\fX_\ph$. To any $X^*\in\bar\cy_\ph$ (see 1.6) and any 
$\bN\in E^2_*\gr^*X^*$ we associate an element $S\in\ce_\ph$ as follows. For any 
$i\in\fX_\ph$ let $Z_i$ be the subspace of $\gr^0X^*$ given by the image of the imbedding 
$\gr^{-a/2}X^*@>>>\gr^0X^*$ induced by $\bN^{a/2}$ for some/any $a\in\ZZ''_{\ge0}$ such 
that $i=f_a$; the natural symplectic form $\la,\ra$ on $\gr^0X^*$ restricts to a symplectic
form on $Z_i$ with $1$-dimensional radical denoted by $L_i$. By definition,
$S=\{(i,j)\in\fX_\ph\T\fX_\ph;L_i=L_j\}$. For $X^*\in\bar\cy_\ph$ and $S\in\ce_\ph$ let
$E^2_*\gr^*X^*_S$ be the set of all $\bN\in E^2_*\gr^*X^*$ such that $(X^*,\bN)$ give rise
to $S$ as above. Let $E^{\ge2}X^*_S=\Ph\i(E^2_*\gr^*X^*_S)$ where $\Ph$ is as in 1.5. We 
thus obtain a partition

$E^{\ge2}_*X^*=\sqc_{S\in\ce_\ph}E^{\ge2}_*X^*_S$
\nl
into finitely many locally closed subvarieties. For $S\in\ce_\ph$ let 

$\Xi_\ph^S=\{g\in SO_Q;g\text{ unipotent }g-1\in E^{\ge2}_*X^*_S\}$.
\nl
Hence $\Xi_\ph$ (see 1.7) is partitioned as $\Xi_\ph=\sqc_{S\in\ce_\ph}\Xi_\ph^S$. Note
that each $\Xi_\ph^S$ is stable under conjugation by $SO_Q$. One can show that the sets 
$\Xi_\ph^S$ with $\ph,S$ as above (together with the sets $\Xi_\ph$ with 
$\ph\in\cf^D-\tcf_1^D$ in the case where $D\in\ZZ''$) are exactly the unipotent conjugacy
classes in $SO_Q$. A quite different classification of unipotent conjugacy classes in
$SO_Q$ is given in \cite{\WA}.

\subhead 4.2\endsubhead
Let $\cn_Q$ be the set of nilpotent elements $\na\in\End(V)$ such that
$\la x,\na y\ra+\la\na x,y\ra=0$ for all $x,y$ in $V$ and $\la x,\na x\ra=0$ for all 
$x\in V$. We can view $\cn_Q$ as the set of nilpotent elements in the Lie algebra of 
$SO_Q$. Note that $SO_Q$ acts by conjugation on $\cn_Q$.

Let $u$ be a unipotent element in $SO_Q$. We associate to $u$ a $SO_Q$-conjugacy class in 
$\cn_Q$ as follows. Let $N=u-1\in\cm_Q$. Let $X^*=V^*_N$. Then $N\in E^{\ge2}X^*$ and 
$\bN\in E^2\gr X^*$ is defined as in 1.5. Let $\bQ$ be the quadratic form on $\gr X^*$ 
defined in 1.5. We have $\bN\in\cn_{\bQ}$. Note that if $D\in\ZZ''$ then the set of 
connected components of $\cj_Q$ and that of $\cj_{\bQ}$ may be naturally identified. We can
find an isomorphism of vector spaces $\gr X^*@>\si>>V$ which carries $\bQ$ to $Q$ and (when
$D\in\ZZ''$) induces the identity map from the set of connected components of $\cj_{\bQ}$ 
to that of $\cj_Q$. This isomorphism carries $\bN$ to an element $\na\in\cn_Q$ whose 
$SO_Q$-orbit is independent of the choice of isomorphism. Note that $1+N\m\na$ defines a 
map
$$\align&\{SO_Q-\text{conjugacy classes of unipotent elements in }SO_Q\}\\&@>>>
\{SO_Q-\text{conjugacy classes of nilpotent elements in }\Lie SO_Q\}.\endalign$$
One can show that this map is injective; it is not in general surjective.

\subhead 4.3\endsubhead
The map 4.2(a) makes sense also in the more general framework of \cite{\LU}. 
Assume that $p>1$. Let $G$ be as in \cite{\LU, 0.1}. We assume that property $\fP_1$ in 
\cite{\LU, 1.1} holds for $G$. Let $u$ be a unipotent element in $G$. By $\fP_1$ we can 
find a unique sequence $\vtr=(G^\vtr_0\sps G^\vtr_1\sps G^\vtr_2\sps\do)$ in $D_G$ such 
that $u\in X^\vtr$ (notation of \cite{\LU, 1.1}). In particular we have $u\in G^\vtr_2$. 
Let $\fg=\Lie G,\fg^\vtr_n=\Lie G^\vtr_n$. Let $\bar u$ be the image of $u$ in 
$G^\vtr_2/G^\vtr_3$. Since $G^\vtr_2/G^\vtr_3$ is a connected commutative unipotent group
in which the $p$-th power of any element is $1$ we see that it is canonically isomorphic to
its Lie algebra $\fg^\vtr_2/\fg^\vtr_3$. Hence $\bar u$ can be identified with an element 
of $\fg^\vtr_2/\fg^\vtr_3$. By definition there exists a homomorphism of algebraic groups
$h:\kk^*@>>>G$ such that if we denote by $\fg^n$ the $n$-eigenspace of the action
$x\m\Ad(h(x))$ of $\kk^*$ on $\fg$ ($n\in\ZZ$) we have $\fg^\vtr_n=\fg^n\op\fg^{n+1}\op\do$
for any $n\ge0$; moreover $h$ is unique up to $G^\vtr_0$-conjugacy. Using the decomposition
$\fg^\vtr_2=\fg^2\op\fg^\vtr_3$ we can identify $\fg^\vtr_2/\fg^\vtr_3$ with $\fg^2$ and we
can view $\bar u$ as an element of $\fg^2$ hence as a nilpotent element of $\fg$. This 
element is well defined up to $G^\vtr_0$-conjugacy. This defines a map

$\{\text{unipotent $G$-conjugacy classes in }G\}@>>>
\{\text{nilpotent $G$-conjugacy classes in }\fg\}$.

\enddocument